\documentclass{siamart190516}
\pdfoutput=1
\usepackage{graphicx}
\usepackage{booktabs}
\usepackage{amsmath}
\usepackage[english]{babel}
\usepackage[utf8]{inputenc}
\usepackage{amssymb}
\usepackage{multicol}
\usepackage[utf8]{inputenc}
\usepackage[english]{babel}
\usepackage{multirow}
\usepackage{hyperref}
\usepackage{cancel} 
\usepackage{bbm}
\newcommand{\R}{\mathbb{R}}

\usepackage{ifthen}
\usepackage{float}
\floatstyle{ruled}
\newfloat{algorithm}{ht}{algo}
\floatname{algorithm}{Algorithm}
\newcounter{algo@row}
\newcounter{algo@rowindent}
\newcommand{\algofont}[1]{\textbf{#1}}
\newcommand{\algonumbersize}[1]{\scriptsize{#1}}
\newcommand{\algopreitem}[1][\arabic{algo@row}]{\texttt{\algonumbersize{#1}}}
\newcommand{\algoitemskip}{\hspace{\value{algo@rowindent}cc}}
\newenvironment{algo}{\vskip.3em\small%
	\begin{list}{\algopreitem\texttt{\algonumbersize{:}}}{%
			\usecounter{algo@row}%
			\setcounter{algo@rowindent}{0}%
			\setlength{\itemindent}{2em}%
			\setlength{\labelwidth}{2em}
			\setlength{\parsep}{0cm}%
		}%
	}{
	\end{list}
}
\newcommand{\algonewnestedopen}[2]{
	\newcommand{#1}[1][]{%
		\ifthenelse{\equal{##1}{}}{\item}{\item[{\algopreitem[##1]}]}
		\algoitemskip\algofont{#2}%
		\addtocounter{algo@rowindent}{1}%
		\ignorespaces
	}
}
\newcommand{\algonewnestedaux}[2]{
	\newcommand{#1}[1][]{
		\addtocounter{algo@rowindent}{-1}
		\ifthenelse{\equal{##1}{}}{\item}{\item[{\algopreitem[##1]}]}
		\algoitemskip\algofont{#2}%
		\addtocounter{algo@rowindent}{+1}%
		\ignorespaces
	}
}
\newcommand{\algonewnestedclose}[2]{
	\newcommand{#1}[1][]{
		\addtocounter{algo@rowindent}{-1}
		\ifthenelse{\equal{##1}{}}{\item}{\item[{\algopreitem[##1]}]}
		\algoitemskip\algofont{#2}%
		\ignorespaces
	}
}
\newcommand{\algonewcommand}[2]{
	\newcommand{#1}[1][default]{
		\ifthenelse{\equal{##1}{default}}{\item}{\item[{\algopreitem[##1]}]}%
		\algoitemskip\algofont{#2}%
		\ignorespaces
	}%
}
\newcommand{\algonewkeyword}[2]{\newcommand{#1}{\algofont{#2}}}
\algonewcommand{\STATE}{\ignorespaces}
\algonewcommand{\INPUT}{Input: }
\algonewcommand{\pINPUT}{\phantom{Input: }}
\algonewcommand{\COMPUTE}{Compute: }
\algonewcommand{\OUTPUT}{Output: }
\algonewcommand{\pOUTPUT}{\phantom{Output: }}
\DeclareMathOperator*{\argmin}{arg\,min}				   
\renewcommand{\t} {^{\top}}								   

\newcommand{\bSigma}{{\boldsymbol{\Sigma}}}

\newcommand{\bepsilon}{{\boldsymbol{\epsilon}}}

\newcommand{\bA}{{\bf A}}
\newcommand{\bB}{{\bf B}}
\newcommand{\bC}{{\bf C}}
\newcommand{\bD}{{\bf D}}

\newcommand{\bF}{{\bf F}}
\newcommand{\bG}{{\bf G}}

\newcommand{\bI}{{\bf I}}
\newcommand{\bJ}{{\bf J}}

\newcommand{\bL}{{\bf L}}

\newcommand{\bP}{{\bf P}}
\newcommand{\bQ}{{\bf Q}}
\newcommand{\bR}{{\bf R}}
\newcommand{\bS}{{\bf S}}

\newcommand{\bU}{{\bf U}}
\newcommand{\bV}{{\bf V}}
\newcommand{\bW}{{\bf W}}
\newcommand{\bX}{{\bf X}}
\newcommand{\bY}{{\bf Y}}
\newcommand{\bZ}{{\bf Z}}

\newcommand{\bb}{{\bf b}}

\newcommand{\bd}{{\bf d}}

\newcommand{\bm}{{\bf m}}

\newcommand{\br}{{\bf r}}
\newcommand{\bs}{{\bf s}}

\newcommand{\bu}{{\bf u}}
\newcommand{\bv}{{\bf v}}
\newcommand{\bw}{{\bf w}}
\newcommand{\bx}{{\bf x}}
\newcommand{\by}{{\bf y}}
\newcommand{\bz}{{\bf z}}

\algonewnestedopen{\IF}{if }
\algonewnestedaux{\ELSEIF}{else if }
\algonewnestedaux{\ELSE}{else }
\algonewnestedclose{\ENDIF}{end if }
\algonewnestedopen{\FOR}{for }
\algonewnestedclose{\ENDFOR}{end for }
\algonewnestedopen{\WHILE}{while }
\algonewnestedclose{\ENDWHILE}{end while }
\algonewcommand{\BREAK}{break}%

\algonewkeyword{\For}{for }%
\algonewkeyword{\To}{to }%
\algonewkeyword{\Do}{do }%
\algonewkeyword{\If}{if }%
\algonewkeyword{\Then}{then }%
\algonewkeyword{\Else}{else }%
\algonewkeyword{\End}{end }%
\algonewkeyword{\AND}{and }%
\algonewkeyword{\True}{true }%
\algonewkeyword{\False}{false }%
\algonewkeyword{\Call}{call }%
\algonewkeyword{\irbleigs}{irbleigs }%
\algonewkeyword{\tridiag}{tridiag}%
\algonewkeyword{\reorth}{reorth}%
\newcommand{\TheTitle}{An $\ell_p$ Variable Projection Method for Large-Scale Separable Nonlinear Inverse Problems} 
\newcommand{\ShortTitle}{$\ell_p$ VarPro for separable nonlinear problems} 
\newcommand{\TheAuthors}{M. Espa\~ nol and M. Pasha}
\headers{\ShortTitle}{\TheAuthors}
\title{{\TheTitle}}
\author{
	Malena I. Espa\~ nol\thanks{School of Mathematical and Statistical Sciences, Arizona State University, Tempe, AZ 85281
		(\email{malena.espanol@asu.edu, mpasha3@asu.edu}).}
	\and
	Mirjeta Pasha$^*${
	}
}

\begin{document}
	\maketitle
	
	\begin{abstract} 
		The variable projection (VarPro) method is an efficient method to solve separable nonlinear least squares problems. In this paper, we propose a modified VarPro for large-scale separable nonlinear inverse problems that promotes edge-preserving and sparsity properties on the desired solution and enhances the convergence of the parameters that define the forward problem. We adopt a majorization minimization method that relies on constructing a quadratic tangent majorant to approximate a general $\ell_p$ regularized problem by an $\ell_2$ regularized problem that can be solved by the aid of generalized Krylov subspace methods at a relatively low cost compared to the original unprojected problem. In addition, we can use more potential general regularizers including total variation (TV), framelet, and wavelets operators. The regularization parameter can be defined automatically at each iteration by means of generalized cross validation. Numerical examples on large-scale two-dimensional imaging problems arising from blind deconvolution are used to highlight the performance of the proposed approach in both quality of the reconstructed image as well as the reconstructed forward operator.
	\end{abstract}
	

	\section{Introduction}
	In many applications in science and engineering ranging from medical to astronomical imaging, the available measurements are contaminated with blur and noise \cite{demmel1997applied,engl1996regularization, deblurring}. We can formulate such problems as discrete ill-posed inverse problems of the form 
	\begin{equation}\label{eq: nonlinear}
		\bG(\by)\bx  = \bd =  \bd_{\rm true} + \bepsilon \quad \mbox{ with } \bG(\by) \bx_{\rm true} = \bd_{\rm true},
	\end{equation}
	where the vector $\bd_{\rm true} \in \mathbb{R}^m$ denotes the unknown error-free vector associated with the available data, and $\bepsilon \in \mathbb{R}^m$ is an unknown error that may stem from measurement or discretization inaccuracies. The matrix $\bG(\by) \in \mathbb{R}^{m \times n}$ models the forward operator that typically is severely ill-conditioned, i.e., it may be numerically rank deficient or may have singular values that cluster at the origin. In this paper, we assume that $\bG$ is unknown, but can be parametrized by a vector $\by\in \mathbb{R}^r$ with $r \ll m$. We call such problems separable nonlinear inverse problems since the observations depend nonlinearly on the vector of unknown parameters $\by$ and linearly on the desired solution $\bx$. Specifically, we aim to compute good approximations of both $\bx$ and $\by$, given a set of observations $\bd = [\bd_1,\bd_2,..., \bd_n]^\top$, a matrix function $\bG$ that maps the unknown vector $\by$ to an $m \times n$ matrix, and an initial approximation $\bx^{(0)}$ of the desired solution $\bx$ as well as an initial approximation $\by^{(0)}$ of the parameter vector $\by$. 
	This leads to a minimization problem of the form
	\begin{equation}\label{eq: min}
		\min_{\bx, \by} \left\|\bG(\by)\bx - \bd\right\|_2^2.
	\end{equation}
	
	Classical optimization approaches that can be used to solve \eqref{eq: min} include block coordinate
	descent, which alternates between fixing one set of variables
	and minimizing with respect to the other set, see for instance \cite{bonettini2011inexact, tseng2009coordinate}. This approach requires a good initial
	guess and is known to have uncertain practical convergence properties. Other methods such as the Gauss-Newton (GN) and Levenberg-Marquardt methods can be used~\cite{nocedal2006numerical}. Among a large set of methods available to solve \eqref{eq: min}, we focus our attention to the variable projection (VarPro) method introduced in \cite{golub1973differentiation}. VarPro is an efficient method where the main idea is to eliminate the linear variable $\bx$ by solving
	a linear least squares (LS) problem for each nonlinear variable $\by$. Therefore, by writing $\bx = \bx(\by) = \bG(\by)^{\dagger}\bd$, where $\bG(\by)^{\dagger}= (\bG(\by)^\top\bG(\by))^{-1}\bG^\top(\by)$ is the Moore-Penrose pseudoinverse of $\bG(\by)$, the functional to be minimized is reduced to a functional of the variable $\by$ only, leading to the following minimization problem
	\begin{equation}\label{eq: minProjec}
		\min_{\by} \|\bG(\by)\bG(\by)^\dagger\bd - \bd\|_2^2.
	\end{equation}
	which is shown in \cite{golub1973differentiation} to have the same solution as \eqref{eq: min}.
	This \emph{reduced} minimization problem can then be solved using GN. For every $\by$, we define the linear operator $\mathcal{P}^{\perp}_{\bG(\by)} = \bI - \bG(\by)\bG(\by)^{\dagger}$ which is the projector onto the orthogonal complement of the column space of $\bG(\by)$. To solve the reduced problem \eqref{eq: minProjec}, an analytic expression of the Jacobian matrix of $\mathcal{P}^{\perp}_{\bG(\by)}\bd$ with respect to the variable $\by$ is given in \cite{golub1973differentiation}, that uses that
	\begin{equation}\label{eq: JacobianGolubPereyra}
		\mathcal{D}{\mathcal{P}_{\bG}^{\perp} = -\mathcal{P}_{\bG}^{\perp}\mathcal{D}(\bG)\bG^{\dagger} - (\mathcal{P}_{\bG}^{\perp}\mathcal{D}(\bG)\bG^{\dagger})^{\top}}, 
	\end{equation}
	where $\mathcal{D}$ represents the derivative operator with respect to the variable $\by$. In \cite{kaufman1975variable}, Kaufman suggested that the second term of \eqref{eq: JacobianGolubPereyra} could be discarded to speed up the calculation of the Jacobian matrix if the residual is small. In \cite{kaufman1978method}, VarPro was extended to constrained problems. In an extended literature, it has been shown computationally and theoretically that separating the linear variables $\bx$ from the nonlinear variable $\by$ by the aid of VarPro speeds up the convergence of iterative methods used to solve problem \eqref{eq: min} (see \cite{golub1976differentiation, kaufman1975variable, ruhe1974algorithms} for a more detailed analysis).
	In \cite{golub2003separable}, Golub and Pereyra reviewed the developments and applications of VarPro to solve separable nonlinear least squares problems that have as their underlying model a linear combination of nonlinear functions. For this same family of problems, O'Leary and Rust implemented a more robust implementation of VarPro and computed the Jacobian matrix more efficiently showing that less iterations are required to converge if the full Jacobian matrix is considered~\cite{o2013variable}.
	
	In the last two decades, several algorithms have been developed to solve the regularized formulation of the nonlinear problem \eqref{eq: min}, which is
	\begin{equation}\label{eq: l2}
		\min_{\bx, \by} \|\bG(\by)\bx - \bd\|_2^2+ \lambda\|\bL\bx\|_2^2,
	\end{equation}
	where $\lambda>0$ is called the \emph{regularization parameter} and $\bL\in \mathbb{R}^{q\times n}$ is a \emph{regularization operator} that is problem dependent. Typical choices of $\bL$ include the identity matrix, discrete first or second derivative operators, or some other discrete transform (e.g., discrete wavelet transform or discrete Fourier transform). In~\cite{chung2010efficient}, Chung and Nagy introduced an efficient method that uses hybrid Krylov subspace approach in order to overcome the high computational cost of solving \eqref{eq: l2} for large-scale inverse problems and for the case when $\bL$ is the identity matrix. In \cite{cornelio2014constrained}, Cornelio et al. applied VarPro to solve non-negative constrained regularized inverse problems. Moreover, they presented an analytic formulation of the Jacobian matrix for the non-negative constrained case. However, in the computations they used the approximate Jacobian from the unconstrained problem with the argument that this approximation slightly affects the accuracy but it saves on computational time. In \cite{herring2017lap}, the authors discussed a linearize and project method based on GN after linearizing the residual that requires the evaluation of the Jacobian matrix and the matrix-vector multiplication of the Jacobian with $\bx$ and $\by$. Determining the regularization parameter using the weighted generalized cross-validation (wGCV) method at every VarPro iteration and efficient approximations of the Jacobian were presented in \cite{chen2018regularized}. More recently, a secant variable projection (SVP) method for solving separable nonlinear least squares problem under the framework of the variable projection that employs rank-one updates to estimate the Jacobian matrices efficiently was proposed in \cite{song2020secant}.
	
	Differently from previous contributions, the focus of this work is on the development of an iterative method that provides edge-preserving and sparsity properties in $\bx$. Edge-preserving and sparsity promoting techniques are well-known in the imaging community and have also been applied to separable nonlinear least squares problems. In \cite{bernstein2020sparse}, the separable nonlinear least squares problem is reformulated as a sparse-recovery problem by the $\ell_1$-norm minimization. The idea has been pioneered in geophysics applications and then it has been widely used \cite{santosa1986linear, taylor1979deconvolution}. In \cite{he2005blind}, He et al. formulated a new time-dependent model for blind deconvolution based on a constrained variational model that uses the sum of the total variation norms as a regularizing functional and both \cite{he2005blind} and \cite{vogel1998fast} addressed the issue of the lack of a unique solution of the blind deconvolution problem.
	
	\subsection{Main contributions of this work} We present an efficient iterative \linebreak method to solve the general separable nonlinear $\ell_2-\ell_p$ minimization problem
	\begin{equation}\label{eq: lp}
		\min_{\bx, \by} \|\bG(\by)\bx - \bd\|_2^2+ \lambda\|\bL\bx\|_p^p.
	\end{equation}
	Notice that the expression 
	\begin{equation*}
		\|\bx\|_p=\left(\sum_{j=1}^n |\bx_{j}|^p\right)^{1/p},\quad \bx=[\bx_1,\bx_2,\ldots,\bx_n]^\top
		\in \mathbb{R}^n
	\end{equation*}
	in the regularization term, is a norm when $p\geq 1$. However, the mapping $\bx\mapsto\|\bx\|_p$ is not a norm for $0<p<1$ since it does not
	satisfy the triangle inequality. Here with a slightly abuse of notation we are going to call $\|\bx\|_p$ a norm for all $p>0$. Because we are interested in reconstructing the edges of $\bx_{\rm true}$, we are particularly interested in values of $0<p\leq 1$ (see Figure \ref{fig: 1D Example}). Even though there is a wide range of methods that solve separable nonlinear inverse problems, there are only a few that provide sparse and edge-preserving properties. To the best of our knowledge, this is the first time that an edge-preserving and sparsity promoting approach is being proposed in a general framework for a general $0<p\leq2$ and with a general regularization matrix $\bL$ that may be a TV, framelet, or wavelet, to mention a few general regularization matrices. 
	
	First, we formulate VarPro for a general $\ell_2$-regularized separable nonlinear least squares problem and carefully derive the Jacobian matrix. For the case when the generalized singular value decomposition (GSVD) of the matrix pair $\{\bG, \bL\}$ is available or feasible to compute, we use it to write the Jacobian matrix in a way that only matrix-vector products are needed. For large scale problems, computing the Jacobian matrix may be computationally not convenient. We showed with a one-dimensional example that the Jacobian of the unregularized problem, which we called \emph{reduced Jacobian}, gives similar convergence curves for both $\bx$ and $\by$ as when the full Jacobian is computed. We also explored the convergence produced by what we called the \emph{half Jacobian}, which consists of the approximation proposed by Kaufman, that discards one of the terms in the full Jacobian under the assumption that the residual is small. We concluded that if the approximate solution is of high quality, no significant difference is noticed in the reconstructions when the reduced and the full Jacobians are used. Therefore, the reduced Jacobian may be preferable specially for large-scale problems. Lastly, we extend VarPro for the case when a closed form solution for $\bx$ is not available which is the case when $p\neq 2$, and show that for large-scale inverse problems, in the step in VarPro where $\by$ is fixed, we solved \eqref{eq: lp} by combining a majorization-minimization strategy to approximate the $\ell_p$ norm by a weighted $\ell_2$. 
	
	This paper is organized as follows. In Section \ref{sec: MM-GKS}, we review the majorize–minimize generalized Krylov subspace (MM-GKS) method to solve the linear problem for the case when $\bG$ is known. In Section \ref{sec: sepNLTik}, we describe how to solve the nonlinear minimization problem for the case when $p=2$ and extend it to $0<p<2$. We describe VarPro for a general $p$ in Section \ref{sec: varpro_lp} by carefully deriving the Jacobian formulations for both cases. Numerical examples illustrate the performance of the proposed methods in Section \ref{sec: numexamples} and some conclusions, future work, and remarks are given in Section \ref{sec: conclusions}.
	
	\section{Iterative method to solve large-scale $\ell_{p}$-regularized linear problems}\label{sec: MM-GKS} Here we discuss an iterative method to solve problems of the form
	\begin{equation}\label{eq: generalLp}
		\min_{\bx}\mathcal{J}_{\lambda,p}(\bx)=\min_{\bx} \|\bG\bx - \bd\|_2^2+ \lambda\|\bL\bx\|_p^p,
	\end{equation}
	for the case when $\bG$ (or its corresponding parametrization vector $\by$) is considered known. There are a huge variety of optimization methods that can be used to solve \eqref{eq: generalLp} when $p=1$, see~\cite{beck2009fast, goldstein2009split, wright2009sparse}. More recent developed iterative methods to solve \eqref{eq: generalLp} for general $0<p\leq 2$ in the case of large-scale inverse problems include iteratively reweighting methods that are based on generalized Krylov subspace methods~\cite{buccini2020modulus, huang2017majorization, lanza2015generalized, rodriguez2008efficient} and on flexible Krylov subspace methods \cite{chung2019flexible, gazzola2014generalized, gazzola2020iteratively}.  Flexible Krylov subspace methods need the computation of the inverse of $\bL$ and that may become computationally demanding to compute for some particular choices of $\bL$. Because we want to consider very general regularization matrices $\bL$, we focus our attention in the group of methods that are based on generalized Krylov subspace methods to iteratively reweight the regularization term in \eqref{eq: generalLp}. More specifically, we choose to apply the majorize–minimize generalized Krylov subspace (MM-GKS) method introduced in \cite{huang2017majorization} to solve  \eqref{eq: generalLp}. First, we will briefly discuss the majorization strategy to approximate \eqref{eq: generalLp} by a sequence of functionals with $\ell_2$ in the regularization term. Let the function $\phi_{p}\colon\R\rightarrow \R_{+}\cup\{+\infty\}$ be defined by
	$\phi_{p}(t)=|t|^{p}$, $p\in \R$. Then, a smoothed approximation of $\phi_{p}(t)$ can be used to approximate the functional $\mathcal{J}_{\lambda,p}$ by a differentiable functional for $p \in (0,1]$. A 
	popular choice is given by 
	\begin{equation*}\label{eq:smothedPhi}
		\phi_{p,\epsilon}(t)= \left(t^2+\epsilon^2\right)^{p/2}\mbox{~~with~~}
		\left\{ \begin{array}{ll}
			\epsilon >0\mbox{~~for~~}0<p\leq 1,\\
			\epsilon =0\mbox{~~for~~}p>1,
		\end{array}
		\right.
	\end{equation*}
	for a small value of $\epsilon>0$. Then, a smoothed approximation of
	$\|\bx\|_p^p$ for a vector $\bx = [\bx_1, \bx_2, ..., \bx_n]\t\in \R^{n}$ can be obtained by
	$\|\bx\|_p^p \approx \sum_{i=1}^{n} \phi_{p,\epsilon}(\bx_i)$ and therefore, instead of solving the $\ell_p$-regularized least squares problem \eqref{eq: generalLp}, we consider solving 
	\begin{equation}\label{eq: smoothedlq}
		\min_{\bx}\mathcal{J}_{\lambda,p,\epsilon}(\bx)=\min_{\bx} \left\|\bG\bx - \bd\right\|_2^2+
		\lambda\sum_{j=1}^{n}\phi_{p,\epsilon}((\bL \bx)_{j}). 
	\end{equation}
	The MM-GKS method constructs a sequence of iterates $\bx^{(k)}$ that converges to a stationary point of the functional $\mathcal{J}_{\lambda, p,\epsilon}(\bx)$. The main idea behind this approach is that the functional $\mathcal{J}_{\lambda, p,\epsilon}(\bx)$ is majorized at each iteration $k$ by a quadratic functional $\mathcal{Q}_{\bx^{(k)}}(\bx)$, whose minimizer provides the next iterate $\bx^{(k+1)}$. To be more precise, we give the following definition. 
	
	\begin{definition}[\cite{huang2017majorization}]\label{def: 1}
		A functional $\mathcal{Q}_{\bv}\colon\R^{n}\rightarrow \R$ is said to be a 
		quadratic tangent majorant for $\mathcal{J}_{\lambda, p,\epsilon}(\bx)$ at $\bv\in\R^n$ if it satisfies
		\begin{enumerate}
			\item $\mathcal{Q}_{\bv}(\bx)$ is quadratic,
			\item $\mathcal{Q}_{\bv}(\bv)=\mathcal{J}_{\lambda, p,\epsilon}(\bv)$, 
			\item $\bigtriangledown_{\bx}\mathcal{Q}_{\bv}(\bv)=\bigtriangledown_{\bx}\mathcal{J}_{\lambda, p,\epsilon}(\bv)$, and
			\item $\mathcal{Q}_{\bv}(\bx)\geq \mathcal{J}_{\lambda, p,\epsilon}(\bx)$ for all $\bx\in \R^{n}$.
		\end{enumerate}
	\end{definition}
	
	Huang et al. \cite{huang2017majorization} describe two approaches to construct a quadratic tangent majorant for \eqref{eq: smoothedlq}. The two different type of majorants are referred to as adaptive or fixed quadratic majorants. The latter are cheaper to compute, but may give slower convergence. In this paper, we focus on the adaptive quadratic majorant, but all results also 
	hold for fixed quadratic majorants. 
	
	Let $\bx^{(k)}$ be the approximate solution at iteration $k$. We define the vector $\bu^{(k)}=\bL\bx^{(k)}$ and the weights
	$\bw_{p,\epsilon}^{(k)} = \left((\bu^{(k)})^2+\epsilon^2\right)^{p/2-1}$, where all operations in the expressions on the right-hand side, including squaring, are 
	element-wise. The derivation of the weights $\bw_{p,\epsilon}^{(k)}$ is given in detail in \cite{lanza2015generalized}. We define the weighting matrix $\bP_{p,\epsilon}^{(k)} = ({\rm diag}(\bw_{p,\epsilon}^{(k)}))^{1/2}$ and consider the quadratic tangent majorant for the functional $\mathcal{J}_{\lambda,p,\epsilon}(\bx)$ at $\bx^{(k)}$ to be
	\begin{equation}\label{eq: QuadraticMajorantQ_2}
		\begin{array}{rcl}
			\mathcal{Q}_{\bx^{(k)}}(\bx)=\displaystyle{\frac{1}{2}} \|\bG\bx-\bd\|^{2}_{2}
			+\displaystyle{\frac{\lambda\epsilon^{p-2}}{2}} \|\bP_{p,\epsilon}^{(k)}\bL\bx\|^{2}_{2}+c,
		\end{array}
	\end{equation}
	where $c$ is a suitable constant that is independent of $\bx^{(k)}$.
	An approximate solution at iteration $k+1$, $\bx^{(k+1)}$, can be determined as the zero of the gradient of \eqref{eq: QuadraticMajorantQ_2} by solving the normal equation
	\begin{equation}\label{eq: normaleqQuadMajorant}
		\left(\bG\t\bG + \eta\bL\t (\bP_{p,\epsilon}^{(k)})^2\bL\right)\bx = \bG\t\bd,
	\end{equation}
	where $\eta = \lambda\epsilon^{p-2}$. The system \eqref{eq: normaleqQuadMajorant} has a unique solution if the condition \linebreak $\mathcal{N}(\bG\t\bG)\cap \mathcal{N}(\bL\t \left(\bP_{p,\epsilon}^{(k)})^2\bL\right)=\{0\}$ is satisfied, which is typically the case. The solution $\bx^{(k+1)}$ of \eqref{eq: normaleqQuadMajorant} is the unique minimizer of the quadratic tangent majorant function $ \mathcal{Q}_{\bx^{(k)}}(\bx)$, hence a solution method is well defined.
	Nevertheless, solving \eqref{eq: normaleqQuadMajorant} for large matrices $\bG$ and $\bL$ may be computationally demanding or even prohibitive. Therefore, the generalized Krylov subspace (GKS) algorithm is applied to solve \eqref{eq: normaleqQuadMajorant}~\cite{lampe2012large}. The GKS method first determines an initial reduction of $\bG$ to a small bidiagonal matrix by applying $1\leq\ell\ll\min\{m,n\}$ steps of Golub--Kahan bidiagonalization to $\bG$ with 
	initial vector $\bd$. This gives the decomposition
	\begin{equation}\label{bdiag}
		\bG\bV_{\ell}=\bU_{\ell}\bB_{\ell},
	\end{equation}
	where the matrix $\bV_{\ell}\in\R^{n\times\ell}$ has orthonormal columns that span the Krylov
	subspace $K_\ell(\bG\t\bG,\bG\t\bd)={\rm span}\{\bG\t\bd,(\bG\t\bG)\bG\t\bd,\ldots,(\bG\t\bG)^{\ell-1}\bG\t\bd\}$, the matrix $\bU_{\ell}\in\R^{m\times(\ell+1)}$ has orthonormal columns, and the matrix
	$\bB_{\ell}\in\R^{(\ell+1)\times\ell}$ is lower bidiagonal. 
	In the problems that we focus on in this paper, the subspace generated by GKS is of small dimension, hence it is inexpensive to compute the SVD or the QR factorization of $\bG\bV_{\ell}$. Therefore, we compute the QR factorization of the tall and skinny matrices and obtain $\bG\bV_{\ell} = \bQ_{\bG}\bR_{\bG}$ and $\bP_{p,\epsilon}^{(0)}\bL\bV_{\ell} = \bQ_{\bL}\bR_{\bL}$, where $\bQ_{\bG}\in\R^{m\times\ell}$ and $\bQ_{\bL}\in\R^{s\times\ell}$ have orthonormal columns and 
	$\bR_{\bG}$, $\bR_{\bL}\in\R^{\ell\times\ell}$ are upper triangular matrices. Here we assume that 
	$1\leq\ell\leq s$ is small enough so that the decomposition $\bG\bV_{\ell} = \bU_{\ell}\bB_{\ell+1, \ell}$
	exists. The corresponding projected least squares problem of \eqref{eq: normaleqQuadMajorant} for $k=1$ is 
	\begin{equation}\label{eq: minKryov2}
		\bz^{(1)}=\arg\min_{\bz\in\R^\ell}\left\|\left[\begin{array}{c} \bR_{\bG}\\ \eta^{1/2}\bR_{\bL}\end{array}\right]
		\bz-\left[\begin{array}{c} \bQ_{\bG}\t\bd\\ 
			0\end{array}\right]\right\|_2^2,
	\end{equation}
	where $\bx^{(1)} = \bV_{\ell}\bz^{(1)}$.
	The normal equation corresponding to \eqref{eq: minKryov2} can be written as
	\begin{equation}\label{eq: normaleqKryl}
		(\bR_{\bG}\t\bR_{\bG} + \eta \bR_{\bL}\t\bR_{\bL})\bz = \bR_{\bG}\t\bQ_{\bG}\bd.
	\end{equation}
	Once $\bz^{(1)}$ is computed, we can compute the residual vector
	\begin{equation*}\label{eq: residual}
		\br^{(1)}=\bG\t(\bG\bV_{\ell} \bz^{(1)} -\bd)+\eta \bL\t\bP_{p,\epsilon}^{(0)}\bL\bV_{\ell}\bz^{(1)}.
	\end{equation*}
	We use the scaled residual vector $\bv_{\rm new}=\br^{(1)}/\|\br^{(1)}\|_2$ to expand the solution subspace. We define the matrix 
	$\bV_{\ell+1}=[\bV_{\ell},\bv_{\rm new}]\in\R^{n\times(\ell+1)}$, whose columns form an orthonormal basis 
	for the expanded solution subspace. Note that in exact arithmetic  $\bv_{\rm new}$ is orthogonal to the columns of $\bV_{\ell}$, but in computer arithmetic the vectors may lose orthogonality, hence reorthogonalization may be  needed. We store the matrices
	\begin{equation*}
		\bG\bV_{\ell+1}=[\bG\bV_{\ell}, \bG\bv_{\rm new}],\quad \bL\bV_{\ell+1}=[\bL\bV_{\ell}, \bL\bv_{\rm new}],
	\end{equation*}
	and iterate obtaining at the $k$th iteration the matrices $\bG \bV_{\ell+k+1} = [\bG \bV_{\ell+k}, \bG\bv_{\rm new}]$ and $\bL\bV_{\ell+k+1} = [\bL\bV_{\ell+k}, \bL\bv_{\rm new}]$. The MM-GKS algorithm to solve \eqref{eq: generalLp} is summarize in Algorithm \ref{Alg: MM-GKS}.
	\begin{algorithm}[!ht]
		\begin{algo}
			\INPUT $\bG, \bL, \bd, \bx^{(0)}, \epsilon$
			\STATE	$\eta=\lambda \epsilon^{p-2}$\;
			\STATE	Generate the initial subspace basis $\bV_{\ell}\in \R^{n\times \ell}$ such that $\bV_{\ell}\t\bV_{\ell}=\bI$\;
			\FOR {$k=0,1,2,\ldots$}{
				\STATE	 $\bu^{(k)}=\bL \bx^{(k)}$\;
				\STATE	 $\bw^{(k)}_{p,\epsilon}=\left((\bu^{(k)})^2+\epsilon^2\right)^{p/2-1}$\;	
				\STATE Compute $\bP_{p,\epsilon}^{(k)} = ({\rm diag}(\bw_{p,\epsilon}^{(k)}))^{1/2}$
				\STATE	Compute $\bG \bV_{\ell+k}$ and $\bP_{p,\epsilon}^{(k)}\bL \bV_{\ell+k}$\;
				\STATE	Compute/Update the QR of $\bG \bV_{\ell+k} = \bQ_{\bG}\bR_{\bG}$ and $\bP_{p,\epsilon}^{(k)}\bL \bV_{\ell+k}=\bQ_{\bL}\bR_{\bL}$\;
				\STATE $\bz^{(k+1)}=( \bR_{\bG}\t\bR_{\bG} + \eta \bR_{\bL}\t\bR_{L})^{-1}\bR_{\bG}\t\bQ_{\bG}\bd $\;
				\STATE $\bx^{(k+1)}=\bV_{\ell+k}\bz^{(k+1)}$\;
				\STATE 
				$
				\br^{(k+1)}=\bG\t(\bG\bV_{\ell+k} \bz^{(k+1)} -\bd)+\eta \bL\t\bP_{p,\epsilon}^{(k)}\bL\bV_{\ell+k}\bz^{(k+1)}
				$\;
				\STATE		Reorthogonalize, if needed, $\br^{(k+1)} = \br^{(k+1)} - \bV_{\ell+k}\bV_{\ell+k}\t\br^{(k+1)}$\;
				\STATE		Enlarge the solution subspace with $\bv_{\rm new}=\frac{\br^{(k)}}{\|\br^{(k)}\|_{2}}$: $\bV_{\ell+k+1}=[\bV_{\ell+k}, \bv_{\rm new}]$\;
			}
			\ENDFOR
		\end{algo}
		\caption{MM-GKS \cite{huang2017majorization}}
		\label{Alg: MM-GKS}
	\end{algorithm}
	
	\subsection{Regularization parameter}\label{sec: regparam}
	The regularization parameter $\lambda$ plays a fundamental role on estimating an approximation of $\bx$ and the right choice of it is very important. Investigations on the importance of the right choice of the regularization parameters are performed in \cite{he2005blind}. The authors in \cite{chung2010efficient} illustrated the importance of choosing a good regularization parameter for each  VarPro iteration applied to solving blind deconvolution and showed how the regularization parameter affects the overall convergence. Different approaches exist to guide the selection of $\lambda$ such as L-curve, discrepancy principle (DP), unbiased predictive risk estimator (UPRE), and generalized cross validation (GCV) \cite{ engl1996regularization, hanke1993regularization, hansen1998rank, vogel2002computational}. Each of them has advantages and disadvantages especially when dealing with large scale problems. For instance, UPRE and DP require a knowledge of the bound of the noise in the available data. Therefore, we will apply the GCV method assuming no knowledge about the level of noise is available.
	
	In \cite{chung2008weighted}, a weighted GCV to compute regularization parameters for Tikhonov regularization was introduced and the SVD of the forward operator $\bG$ was used to compute the GCV function \eqref{eq: GCV} efficiently. However, computing the SVD might not be computationally attractive or might even prohibitive for large scale problems, and therefore, they also computed the GCV function by the aid of a Lanczos type method. Other references where efficient ways to compute the GCV can be found in \cite{buccini2021generalized, fenu2016gcv, fenu2017gcv}. 
	For problems of the form like \eqref{eq: generalLp} with $p=2$, GCV chooses the value $\lambda$ that minimizes the function 
	\begin{equation}\label{eq: GCV}
		\mathcal{G}_{\bG, \bd}(\lambda) = \frac{n \left\|(\bI - \bG\bG_\bL^{\dagger}\bd \right\|^2_2}{\left(\rm trace \left(\bI - \omega\bG\bG_\bL^{\dagger}\right)\right)^2},
	\end{equation}
	where the weighting factor $\omega$, $0<\omega\leq 1$ can be modified so that larger values of $\omega$ result in larger regularization parameters, and smaller values of $\omega$ result in smaller values of $\lambda$ computed.
	
	Here, we need to apply GCV to find $\eta$ (and therefore $\lambda$) in each iteration of MM-GKS. Then, the GCV function \eqref{eq: GCV} applied to \eqref{eq: minKryov2} can be written as
	\begin{equation}\label{eq: GCVKryl}
		\mathcal{G}_{\bR_{\bG}, \hat{\bd}}(\eta) = \frac{(\ell+k) \left\|(\bI - \bR_{\bG}\bR_{\bG_{\eta}}^{\dagger})\bQ_\bG\t\bd \right\|^2_2}{\left(\rm trace \left(\bI - \omega \bR_{\bG}\bR_{\bG_{\eta}}^{\dagger}\right)\right)^2},
	\end{equation}
	where $\bR_{\bG_\eta}^{\dagger}=(\bR_{\bG}\t\bR_{\bG} + \eta \bR_{\bL}\t\bR_{\bL})\bR_{\bG}\t$.
	
	Since the matrices $\bR_{\bG}$ and $\bR_{\bL}$ are not very large, we then compute the GSVD of them, which is define by
	\begin{equation}\label{eq: GSVD}
		\bR_{\bG} = \bX_{\bG}\bSigma_{\bG}\bY\t, \qquad
		\bR_{\bL} = \bX_{\bL}\bSigma_{\bL}\bY\t,
	\end{equation}
	where the matrices $\bX_{\bG}\in\R^{(\ell+k)\times (\ell+k)}$ and $\bX_{\bL}\in\R^{(\ell+k)\times (\ell+k)}$ are orthogonal, the matrix 
	$\bY\in\R^{(\ell+k)\times (\ell+k)}$ is nonsingular, and 
	\begin{small}\begin{eqnarray*}
			\bSigma_{\bG}={\rm diag}(\sigma_1,\sigma_2,\ldots,\sigma_{\ell+k})\in\R^{(\ell+k) \times (\ell+k)},\quad
			\bSigma_{\bL}={\rm diag}(\lambda_1,\lambda_2,\ldots,\lambda_{\ell+k})\in\R^{(\ell+k) \times (\ell+k)}.
		\end{eqnarray*}
	\end{small}
	Substituting \eqref{eq: GSVD} into \eqref{eq: GCVKryl} we obtain the regularization parameter $\eta$ by solving the following minimization problem
	\begin{equation}\label{eq: GCVgsvd}
		\eta^{(k+1)}= \arg\min_{\eta} \frac{\left\|\bX_{\bG}\t\hat{\bd}-\bSigma_{\bG}\t(\bSigma_{\bG}\t\bSigma_{\bG} + \eta\bSigma_{\bL}\t\bSigma_\bL)^{-1}\bSigma_{\bG}\t\bX_{\bG}\t\hat{\bd}\right\|_2^2}{\left(\rm trace \left(\bI - \omega\bSigma\t(\bSigma_{\bG}\t\bSigma_{\bG} + \lambda \bSigma_{\bL}\t\bSigma_{\bL})^{-1}\bSigma_{\bG} \right)\right)^2}.
	\end{equation}
	We will not dwell in more details on the regularization parameter as we direct the reader to \cite{buccini2021generalized}.
	
	\section{Iterative methods for $\ell_2$-regularized nonlinear least squares problems}\label{sec: sepNLTik}
	In this section we focus our attention to solve \eqref{eq: lp} for $p=2$ and a general regularization matrix $\bL$. First, we treat \eqref{eq: lp} as a non separable nonlinear least squares problem and apply the Gauss-Newton (GN) method~\cite{dennis1996numerical, kelley1999iterative, nocedal2006numerical} to solve  
	\begin{equation}\label{eq: sepNLTik}
		\min_{\bx, \by}\left\|\left[\begin{array}{c} \bG(\by) \\ \lambda^{1/2} \bL\end{array}\right] \bx-
		\left[\begin{array}{c} \bd \\ \bf{0} \end{array}\right]\right\|^2_2,
	\end{equation}
	or equivalently
	\begin{equation}\label{eq: minF}
		\min_{\bm}{\mathcal{F}(\bm)} = \min_{\bm}\left\|\bF(\bm)\right\|_2^2,
	\end{equation}
	where $\bm=(\bx,\by)$ and $\bF \colon \mathbb{R}^{n+r} \to \mathbb{R}^{m+q}$ is defined by
	$$
	\bF(\bm) = \bF(\bx, \by) =  \left[\begin{array}{c} \bG(\by) \\ \lambda^{1/2} \bL\end{array}\right] \bx-
	\left[\begin{array}{c} \bd \\ \bf{0} \end{array}\right].
	$$
	Thus, GN iterations are defined by
	$$
	\bm^{(i+1)} = \bm^{(i)} + \bs^{(i)}, \,i = 1,2,...,
	$$ 
	where $\bm^{(i)}=\left(\bx^{(i)},\by^{(i)}\right)$ and $\bs^{(i)}$ is the solution of 
	\begin{equation}\label{eq: GN}
		\bs^{(i)} = \argmin_{\bs}\left\| \bJ_{\bF}(\bm^{(i)})\bs - \bF(\bm^{(i)}) \right\|^2_2
	\end{equation}
	with $\bJ_{\bF}\colon \mathbb{R}^{n+r}\to\mathbb{R}^{m+q}$ being the Jacobian matrix of $\bF$ defined by 
	$$\bJ_{\bF}(\bm)=\bJ_{\bF}(\bx,\by) = \left[ \frac{\partial \bF(\bx,\by)}{\partial \bx} , \frac{\partial \bF(\bx,\by)}{\partial \by}\right]=\begin{bmatrix} \bG(\by) & \frac{\partial}{\partial \by} [\bG(\by) \bx] \\ \lambda^{1/2}\bL & \bf{0}\end{bmatrix}.$$
	Notice that each column of $\frac{\partial}{\partial \by} [\bG(\by) \bx]$ corresponds to the partial with respect to each $\by_j$, $j=1,\dots, q$, which is the column vector defined by $\frac{\partial \bG(\by)}{\partial \by_j} \bx$.
	
	\begin{algorithm}[ht]
		\begin{algo}
			\INPUT $\bm^{(0)} = (\bx^{(0)}, \by^{(0)})$, $\bG(\by^{(0)})$, $\bL$, and $\bd$.
			\FOR $i = 0, 1,\dots $
			\STATE $\bF^{(i)} =
			\left[\begin{array}{c} \bG(\by^{(i)}) \\ \lambda^{1/2} \bL\end{array}\right] \bx^{(i)}-\left[\begin{array}{c} \bd \\ \bf{0} \end{array}\right]$
			\STATE Compute/Update the Jacobian matrix $\bJ^{(i)}_{\bF}$
			\STATE $\bs^{(i)} = \argmin_{\bs}\left\|\bJ^{(i)}_{\bF}\bs -\bF^{(i)}\right\|^2_2$
			\STATE $\bm^{(i+1)} = \bm^{(i)} + \bs^{(i)}$
			\ENDFOR
		\end{algo}
		\caption{GN-NLS Algorithm}
		\label{Alg: G_GaussNewton}
	\end{algorithm} 
	
	In Algorithm \ref{Alg: G_GaussNewton} we summarize GN applied to \eqref{eq: sepNLTik}. The main disadvantages are a) to choose the regularization parameter $\lambda$ which is a task that may be computationally demanding for large-scale problems, and b) the coupling of $\bx$ and $\by$ does not actually help in their convergences as has been shown in the literature before. Therefore, we will follow the idea behind VarPro~\cite{golub1973differentiation} that completely eliminates the variable $\bx$ from the problem formulation and provides us with a reduced functional to minimize only with respect to $\by$. That is to say, we can apply GN method to the functional $\Psi(\by) = \mathcal{F}(\bx(\by), \by)$, where $\bx(\by)$ is the solution of the minimization problem
	\begin{equation}\label{eq: Fx}
		\min_{\bx}\mathcal{F}(\bx, \by) =  \min_{\bx}\left\|\left[\begin{array}{c} \bG(\by) \\ \lambda^{1/2} \bL\end{array}\right] \bx-
		\left[\begin{array}{c} \bd \\ \bf{0} \end{array}\right]\right\|_2^2.
	\end{equation}
	A big advantage of doing so is that the choice of the regularization parameter is shifted to the linear problem and therefore well-known regularization selection methods can be applied. For large-scale inverse problems and $\bL = \bI$, we recall that Chung and Nagy solved \eqref{eq: Fx} by the aid of a hybrid Krylov method~\cite{chung2010efficient}. For the case of a more general regularization matrix $\bL$, we can apply  projection-based approaches designed for general-form Tikhonov~\cite{kilmer2007projection, lampe2012large}. 
	
	Notice that a closed form solution of \eqref{eq: Fx} is available, which can be used to rewrite the nonlinear problem with respect to the parameters $\by$, obtaining the reduced minimization problem
	\begin{equation}\label{eq: Psiy}
		\min_{\by}\Psi(\by) = \min_{\by}\|\tilde \bF(\by)\|_2^2, 
	\end{equation}
	where
	\begin{align}
		\tilde \bF(\by) &=  \left[\begin{array}{c} \bG(\by) \\ \lambda^{1/2} \bL\end{array}\right] \left[\begin{array}{c} \bG(\by) \\ \lambda^{1/2} \bL\end{array}\right]^\dagger\left[\begin{array}{c} \bd \\ \bf{0} \end{array}\right]-
		\left[\begin{array}{c} \bd \\ \bf{0} \end{array}\right] \nonumber \\
		&= \left(\left[\begin{array}{c} \bG(\by) \\ \lambda^{1/2} \bL\end{array}\right] \left[\begin{array}{c} \bG(\by) \\ \lambda^{1/2} \bL\end{array}\right]^\dagger - \bI
		\right) \left[\begin{array}{c} \bd \\ \bf{0} \end{array}\right].
	\end{align}
	To simplify notation we define
	\[ \bG_\bL(\by):=\left[\begin{array}{c} \bG(\by) \\ \lambda^{1/2} \bL\end{array}\right]  \mbox{ and }  \mathcal{P}^\perp_{\bG_\bL}(\by):= \bI-\bG_{\bL}(\by)\bG_{\bL}^\dagger(\by),\]
	and write only $\bG_\bL$ and $\mathcal{P}^\perp_{\bG_\bL}$ instead of $\bG_\bL(\by)$ and $\mathcal{P}^\perp_{\bG_\bL}(\by)$ for even more simplification.
	To solve \eqref{eq: Psiy}, the GN iterations are defined by
	$$
	\by^{(i+1)} = \by^{(i)} + \tilde\bs^{(i)}, \,i = 1,2,...,
	$$ 
	where $\tilde\bs^{(i)}$ is the solution of 
	\begin{equation}\label{eq: bigJacobian}
		\tilde\bs^{(i)} = \argmin_{\tilde\bs}\left\| \bJ_{\tilde\bF}(\by^{(i)})\tilde\bs + \tilde\bF(\by^{(i)}) \right\|^2_2
	\end{equation}
	with $\bJ_{\tilde\bF}\colon \mathbb{R}^{r}\to\mathbb{R}^{m+q}$ being the Jacobian matrix of $\tilde\bF$. The $j$-th column of  $\bJ_{\tilde\bF}$ can be computed by
	
	\begin{align}\label{eq: jthColumnJacobianstart}
		\frac{\partial \tilde\bF(\by)}{\partial \by_j} & = \frac{\partial }{\partial \by_j}\left(\bG_{\bL}\bG_{\bL}^\dagger \left[\begin{array}{c} \bd \\ \bf{0} \end{array}\right]\right) 
		=\frac{\partial }{\partial \by_j}\left(\bG_{\bL}\bG_{\bL}^\dagger\right) \left[\begin{array}{c} \bd \\ \bf{0} \end{array}\right] \\ 
		& = \left(\frac{\partial \bG_{\bL} }{\partial \by_j}\bG_{\bL}^\dagger+ \bG_{\bL}\frac{\partial \bG_{\bL}^\dagger }{\partial \by_j} \right) \left[\begin{array}{c} \bd \\ \bf{0} \end{array}\right].
	\end{align}
	By writing $\bG_{\bL}^\dagger=(\bG_{\bL}^\top\bG_{\bL})^{-1}\bG_{\bL}^\top $, applying the product rule, and using that \linebreak $\frac{\partial \bG_{\bL}^{-1}}{\partial \by_j}= -\bG_{\bL}^{-1} \frac{\partial \bG_{\bL}}{\partial \by_j} \bG_{\bL}^{-1}$, we have that
	\begin{align}
		\frac{\partial \bG_{\bL}^\dagger}{\partial \by_j} & = (\bG_{\bL}^\top\bG_{\bL})^{-1} \frac{\partial \bG_{\bL}^\top}{\partial \by_j}\mathcal{P}^\perp_{\bG_{\bL}}-\bG_{\bL}^\dagger \frac{\partial \bG_{\bL}}{\partial \by_j} \bG_{\bL}^\dagger. \nonumber \\
		& = (\bG_{\bL}^\top\bG_{\bL})^{-1} \left[ \frac{\partial \bG^\top}{\partial \by_j}, \bf{0} \right]\mathcal{P}^\perp_{\bG_\bL} \nonumber - \bG_{\bL}^\dagger \left[\begin{array}{c} \frac{\partial \bG}{\partial \by_j} \\ \bf{0} \end{array} \right] \bG_{\bL}^\dagger.
	\end{align} 
	Therefore, the $j$-th column of the Jacobian is given by
	
	\begin{align}\label{eq: jthColumnJacobian}
		\frac{\partial \tilde\bF(\by)}{\partial \by_j} & = 
		\left(\frac{\partial \bG_{\bL} }{\partial \by_j}\bG_{\bL}^\dagger +  \bG_{\bL}\left((\bG_{\bL}^\top\bG_{\bL})^{-1} \left[ \frac{\partial \bG^\top}{\partial \by_j}, \bf{0} \right]\mathcal{P}^\perp_{\bG_\bL} \nonumber - \bG_{\bL}^\dagger \left[ \begin{array}{c} \frac{\partial \bG}{\partial \by_j} \\ \bf{0} \end{array} \right] \bG_{\bL}^\dagger \right) \right) \left[\begin{array}{c} \bd \\ \bf{0} \end{array}\right]\nonumber \\
		& = - \left( \mathcal{P}^\perp_{\bG_\bL} \left[ \begin{array}{c} \frac{\partial \bG}{\partial \by_j} \\ \bf{0} \end{array} \right] \bG_{\bL}^\dagger +
		\left( \mathcal{P}^\perp_{\bG_\bL} \left[ \begin{array}{c} \frac{\partial \bG}{\partial \by_j} \\ \bf{0} \end{array} \right] \bG_{\bL}^\dagger \right)^\top\right)\left[\begin{array}{c} \bd \\ \bf{0} \end{array}\right] \nonumber\\
		& = - \mathcal{P}^\perp_{\bG_\bL} \left[ \begin{array}{c} \frac{\partial \bG}{\partial \by_j} \\ \bf{0} \end{array} \right] \bG_{\bL}^\dagger \left[\begin{array}{c} \bd \\ \bf{0} \end{array}\right] 
		-
		\left( \mathcal{P}^\perp_{\bG_\bL}\left[ \begin{array}{c} \frac{\partial \bG}{\partial \by_j} \\ \bf{0} \end{array} \right] \bG_{\bL}^\dagger \right)^\top\left[\begin{array}{c} \bd \\ \bf{0} \end{array}\right] \nonumber \\ 
		& = -\bA_j - \bB_j.
	\end{align}
	
	Suppose the GSVD of the matrix pair $\{\bG,\bL\}$ is available. Then, we have that $\bG = \bP\bC\bZ\t$ and $\bL = \bar \bP\bS\bZ\t,$ where the matrices $\bP\in\R^{m\times m}$ and $\bar \bP\in\R^{q\times q}$ are orthogonal, the matrix 
	$\bZ\in\R^{n\times n}$ is nonsingular, and
	$\bC ={\rm diag}(c_1,\dots,c_q,1,1,\ldots,1)\in\R^{m\times n}$ and $\bS={\rm diag}(s_1,\dots,s_q,0,0,\ldots,0)\in\R^{q\times n},$
	where we assume that $1\leq q\leq n$. Using the GSVD, we obtain the following expressions of $\bG_{\bL}^\dagger$ and $\mathcal{P}^\perp_{\bG_{\bL}}$
	\begin{align}
		\bG_{\bL}^\dagger & = \bZ^{-\top}(\bC\t\bC + \lambda \bS\t\bS)^{-1}\begin{bmatrix}\bC\t \bP\t, & \lambda^{1/2} \bS\t \bar \bP\t 
		\end{bmatrix} \\
		\mathcal{P}^\perp_{\bG_{\bL}} &= \bI - \begin{bmatrix}\bP\bC \\ \lambda^{1/2} \bar \bP\bS
		\end{bmatrix} (\bC\t\bC + \lambda \bS\t\bS)^{-1} \begin{bmatrix}\bC\t \bP\t, & \lambda^{1/2} \bS\t \bar \bP\t 
		\end{bmatrix}. 
	\end{align}
	Therefore, each term in the $j$-th column of the Jacobian can be computed as follows.
	\begin{align} -\bA_j & =  - \mathcal{P}^\perp_{\bG_\bL} \left[ \begin{array}{c} \frac{\partial \bG}{\partial \by_j} \\ \bf{0} \end{array} \right] \bx^{(i)} \nonumber \\
		& = -   \left[ \begin{array}{c} \frac{\partial \bG}{\partial \by_j} \\ \bf{0} \end{array} \right] \bx^{(i)} +  \begin{bmatrix}\bP\bC \\ \lambda^{1/2} \bar \bP\bS 
		\end{bmatrix}\begin{bmatrix}\bD_\bC \bP\t, & \lambda^{1/2} \bD_\bS \bar \bP\t 
		\end{bmatrix} \left[ \begin{array}{c} \frac{\partial \bG}{\partial \by_j} \\ \bf{0} \end{array} \right] \bx^{(i)} \nonumber \\
		& = -   \left[ \begin{array}{c} \frac{\partial \bG}{\partial \by_j} \\ \bf{0} \end{array} \right] \bx^{(i)} +  \begin{bmatrix}\bP\bC \\ \lambda^{1/2} \bar \bP\bS 
		\end{bmatrix}\bD_\bC \bP\t \frac{\partial \bG}{\partial \by_j}\bx^{(i)} \nonumber \\
		-\bB_j & =
		- (\bG_{\bL}^\dagger) ^\top \left[ \frac{\partial \bG^\top}{\partial \by_j}, \bf{0} \right] \mathcal{P}^\perp_{\bG_\bL} \left[\begin{array}{c} \bd \\ \bf{0} \end{array}\right]  \nonumber \\
		& =
		\begin{bmatrix}\bP\bD_\bC \\ \lambda^{1/2} \bar \bP\bD_\bS 
		\end{bmatrix}\bZ^{-1} \left[ \frac{\partial \bG^\top}{\partial \by_j}, \bf{0} \right] \bF(\by^{(i-1)}), \nonumber \\
		& =
		\begin{bmatrix}\bP\bD_\bC \\ \lambda^{1/2} \bar \bP\bD_\bS 
		\end{bmatrix}\bZ^{-1} \frac{\partial \bG^\top}{\partial \by_j} (\bG(\by^{(i-1)})\bx^{(i)}-\bd),\nonumber
	\end{align}
	where $\bD_\bC$ and $\bD_\bS$ are the diagonal matrices defined by $\bD_\bC = \bC(\bC\t\bC + \lambda \bS\t\bS)^{-1}$ and $\bD_\bS =\bS (\bC\t\bC + \lambda \bS\t\bS)^{-1}$, and $\bx^{(i)}$ is the solution of the least squares problem with $\by^{(i-1)}$ computed in the previous VarPro iteration. Here we have extended the formulation presented in \cite{o2013variable}, so that only matrix-vector multiplications are used and the only inverse to be computed is one corresponding to a diagonal matrix. In Algorithm \ref{Alg: Gen-VarPro} we summarize VarPro applied to \eqref{eq: sepNLTik} for $p=2$.
	
	\begin{algorithm}[!ht]
		\begin{algo}
			\INPUT $\bG(\by^{(0)})$, $\bL$, and $\bd$.
			\FOR $i = 0, 1,... $
			\STATE $\bx^{(i)} = \min_{\bx} \left\| \left[\begin{array}{c} \bG(\by^{(i)}) \\ \lambda^{1/2} \bL\end{array}\right] \bx-
			\left[\begin{array}{c} \bd \\ \bf{0} \end{array}\right]\right\|_2^2$
			\STATE $\tilde\bF^{(i)} =
			\left[\begin{array}{c} \bG(\by^{(i)}) \\ \lambda^{1/2} \bL\end{array}\right] \bx^{(i)} -\left[\begin{array}{c} \bd \\ \bf{0} \end{array}\right]$
			\STATE Compute/Update the Jacobian matrix $\bJ^{(i)}_{\tilde\bF}$
			\STATE $\tilde\bs^{(i)} = \argmin_{\bs}\left\|\bJ^{(i)}_{\tilde\bF}\tilde\bs - \tilde\bF^{(i)}\right\|^2_2$
			\STATE $\by^{(i+1)} = \by^{(i)} + \tilde\bs^{(i)}$
			\ENDFOR
		\end{algo}
		\caption{Generalized Variable Projection (Gen-VarPro)}
		\label{Alg: Gen-VarPro}
	\end{algorithm} 
	
	\section{$\ell_p$ VarPro for large-scale inverse problems}\label{sec: varpro_lp}
	Now we look at the case when $0<p<2$ and present the $\ell_p$ VarPro method. Recall that the minimization problem that we aim to solve is \begin{equation}\label{eq: sepNLLp}
		\min_{\bx, \by}\left\|\bG(\by)\bx - \bb\right\|^2_2 + \lambda\|\bL\bx\|_p^p.
	\end{equation}
	Differently from the case when $p=2$, a closed form solution of \eqref{eq: sepNLLp} for $\by$ fixed does not exist. Nevertheless, we consider an approximation of the solution $\bx(\by)$.  We define the function $\hat{\Psi}(\by) = \mathcal{\hat{F}}(\bx(\by), \by)$, where $\bx(\by)$ is the solution of the minimization problem
	\begin{equation}\label{eq: hatFx}
		\min_{\bx}\mathcal{\hat{F}}(\bx, \by) =  \min_{\bx}\left\|\bG(\by)\bx - \bd\right\|_2^2 + \lambda\|\bL\bx\|_p^p.
	\end{equation}
	Let $\bx^{(k)}$ be an available approximate solution of \eqref{eq: hatFx}. We define the vector $\bu^{(k)}=\bL\bx^{(k)}$ and the weights
	$\bw_{p,\epsilon}^{(k)} = \left((\bu^{(k)})^2+\epsilon^2\right)^{p/2-1}$, where all operations in the expressions on the right-hand side, including squaring, are 
	element-wise. We define the weighting matrix $\bP_{p,\epsilon}^{(k)} = ({\rm diag}(\bw_{p,\epsilon}^{(k)}))^{1/2}$ and consider the quadratic tangent majorant for the functional $\mathcal{J}_{\lambda,p,\epsilon}(\bx)$ at $\bx^{(k)}$ to be
	\begin{equation}\label{eq: QuadraticMajorantQ}
		\begin{array}{rcl}
			\mathcal{Q}_{\bx^{(k)}}(\bx)=\displaystyle{\frac{1}{2}} \|\bG(\by)\bx-\bd\|^{2}_{2}
			+\displaystyle{\frac{\lambda\epsilon^{p-2}}{2}} \|\underbrace{\bP_{p,\epsilon}^{(k)}\bL}_{ \hat{\bL}}\bx\|^{2}_{2}+c.
		\end{array}
	\end{equation}
	Computing the gradient with respect to $\bx$ results in minimizing the following: 
	\begin{equation}\label{eq: normaleqQuadMajorantCase2}
		\left(\bG(\by)\t\bG(\by) + \eta\bL\t (\bP_{p,\epsilon}^{(k)})^2\bL\right)\bx = \bG(\by)\t\bd,
	\end{equation}
	With this choice of $\bP_{p,\epsilon}^{(k)}$ we seek for an approximation  such that $\|\bL\bx\|^{p}_{p}\approx\|\bP_{p,\epsilon}^{(k)}\bL\bx\|^{2}_{2}$. 
	
	We want to adapt the gen-VarPro for $p\neq 2$. Then, we consider the reduced functional
	\begin{equation}\label{eq: hatPsiy}
		\min_{\by}\hat{\Psi}(\by) = \min_{\by}\|\hat \bF(\by)\|_2^2, 
	\end{equation}
	where
	\begin{align}
		\hat \bF(\by) &=  \left[\begin{array}{c} \bG(\by) \\ \lambda^{1/2}\bP_{p,\epsilon}^{(k)} \bL\end{array}\right] \left[\begin{array}{c} \bG(\by) \\ \lambda^{1/2} \bP_{p,\epsilon}^{(k)}\bL\end{array}\right]^\dagger\left[\begin{array}{c} \bd \\ \bf{0} \end{array}\right]-
		\left[\begin{array}{c} \bd \\ \bf{0} \end{array}\right] \nonumber \\
		&= \left(\left[\begin{array}{c} \bG(\by) \\ \lambda^{1/2}\bP_{p,\epsilon}^{(k)} \bL\end{array}\right] \left[\begin{array}{c} \bG(\by) \\ \lambda^{1/2} \bP_{p,\epsilon}^{(k)} \bL\end{array}\right]^\dagger - \bI
		\right) \left[\begin{array}{c} \bd \\ \bf{0} \end{array}\right].
	\end{align}
	Notice that the only difference is that we have $\bP_{p,\epsilon}^{(k)} \bL$ instead of $\bL$. Our approach is based on the following observations. Computing $\bx^{(i)}$ for large scale problems might be computationally expensive or even prohibitive, hence we apply MM-GKS to simultaneously compute an approximate solution $\bx^{(i)}$ of the desired solution $\bx$ and automatically define the regularization parameter without too much effort. Once $\bx^{(i)}$ is computed, we can define $\bP_{p,\epsilon}^{(i)}$ and therefore redefine the regularization term to $\hat\bL = \bP_{p,\epsilon}^{(i)}\bL$ and compute the Jacobian as in \eqref{eq: jthColumnJacobianstart}-\eqref{eq: jthColumnJacobian}. In Algorithm \ref{Alg: lp VarPro}, we summarize all the steps to solve the $\ell_p$-regularized separable nonlinear problem \eqref{eq: sepNLLp}.
	
	\begin{algorithm}
		\begin{algo}
			\INPUT $\bG, \bL, \bd, \bx^{(0)}$,
			\FOR {$i=0,1,2,\ldots$}{
				\STATE $\bx^{(i)} = \argmin_{\bx}\|\bG\bx - \bd\|_2^2+ \lambda\|\bL\bx\|_p^p,$ (solved using MM-GKS) 
				\STATE	 $\bu^{(i)}=\bL \bx^{(i)}$\;
				\STATE	 $\bw^{(i)}_{p,\epsilon}=\left((\bu^{(i)})^2+\epsilon^2\right)^{p/2-1}$\;	
				\STATE Compute $\bP_{p,\epsilon}^{(i)} = ({\rm diag}(\bw_{p,\epsilon}^{(i)}))^{1/2}$
				\STATE Compute $\hat\bL = \bP_{p,\epsilon}^{(i)}\bL$ 
				\STATE $\hat\bF^{(i)} =
				\left[\begin{array}{c} \bG(\by^{(i)}) \\ \lambda^{1/2} \hat\bL\end{array}\right] \bx^{(i)} -\left[\begin{array}{c} \bd \\ \bf{0} \end{array}\right]$
				\STATE Compute/Update the Jacobian matrix $\bJ^{(i)}_{\hat\bF}$
				\STATE $\tilde\bs^{(i)} = \argmin_{\bs}\left\|\bJ^{(i)}_{\hat\bF}\tilde\bs - \hat\bF^{(i)}\right\|^2_2$
				\STATE $\by^{(i+1)} = \by^{(i)} + \bs^{(i)}$}
			\ENDFOR
		\end{algo}
		\caption{NLS$_{\ell_2-\ell_p}$}
		\label{Alg: lp VarPro}
	\end{algorithm}
	
	\section{Numerical experiments}\label{sec: numexamples}
	In this section, we present a few numerical examples to demonstrate the performance of the 
	methods discussed in the previous sections. We start with a one-dimensional deconvolution problem, which allows us to easily compute and test the different Jacobian matrices. We then show some numerical results for two-dimensional blind deconvolution problems where we test the performance of the $\ell_p$ VarPro with the inclusion of the MM-GKS method to find an approximation of $\bx$. We use the MATLAB package \cite{nagy2004iterative} to efficiently compute the Jacobian approximation. All tests were performed using MATLAB R2020b on a single processor, Intel Core i7 computer. To assess the quality of the reconstructed solution at iteration $i$, we compute the Relative Reconstruction Error (RRE) defined  by
	\begin{equation}\label{RRE}
		{\rm RRE}(\bx^{(i)})=\frac{\|\bx^{(i)}-\bx_{\rm true}\|_2}{\|\bx_{\rm true}\|_2}.
	\end{equation}
	To assess the quality of the reconstructed parameters $\by$ at iteration $i$ we use the RRE as 
	\begin{equation}\label{RREy}
		{\rm RRE}(\by^{(i)})=\frac{\|\by^{(i)}-\by_{\rm true}\|_2}{\|\by_{\rm true}\|_2}.
	\end{equation}
	Visual inspection of the reconstructions, quality measures, as well as the actual parameters $\by$ recovered are also reported.
	\begin{figure}[h!]
		\centering
		\begin{minipage}{0.32\textwidth}
			\centering
			\includegraphics[width = \textwidth]{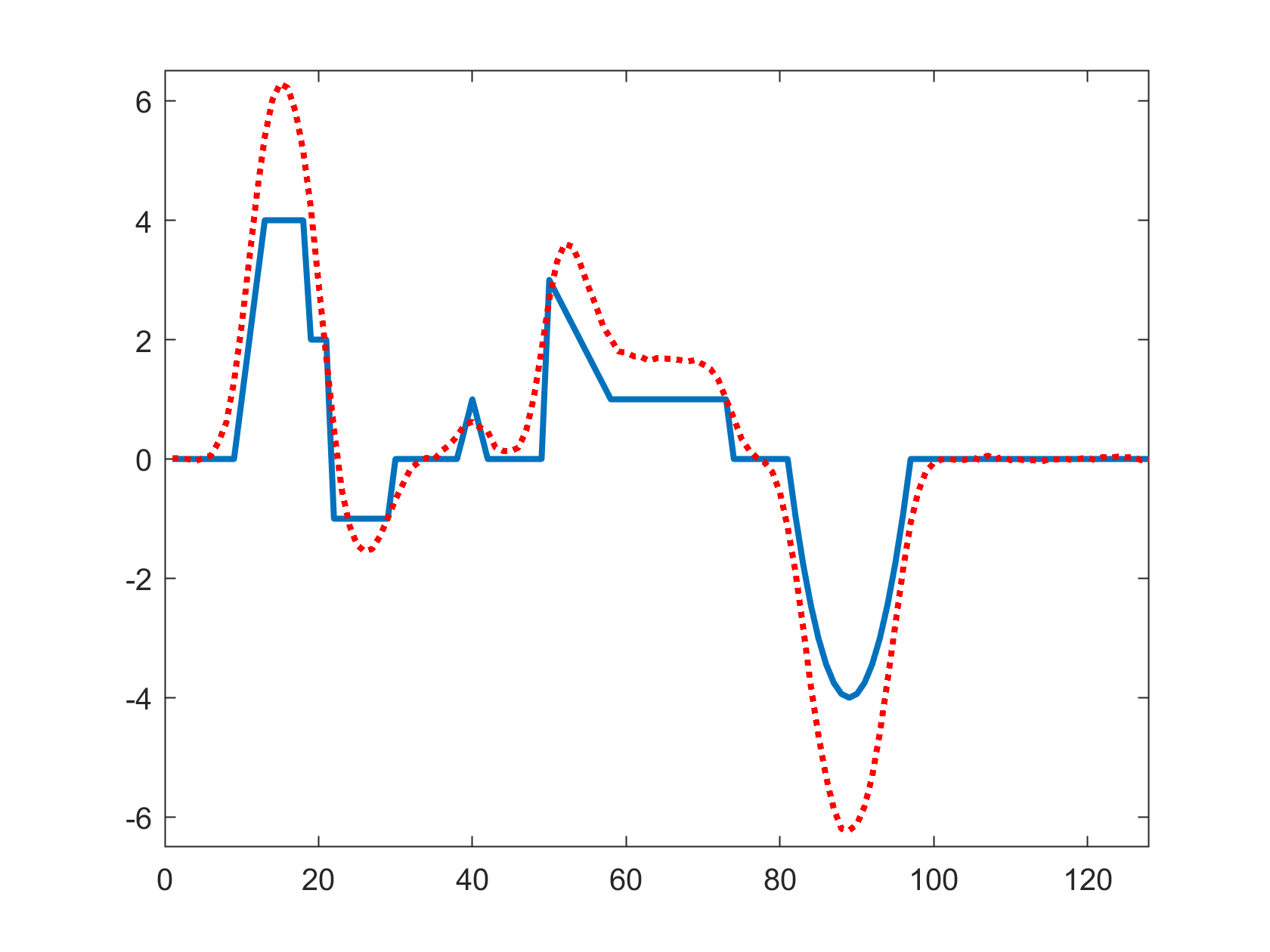}\\a)
		\end{minipage}
		\begin{minipage}{0.32\textwidth}
			\centering
			\includegraphics[width = \textwidth]{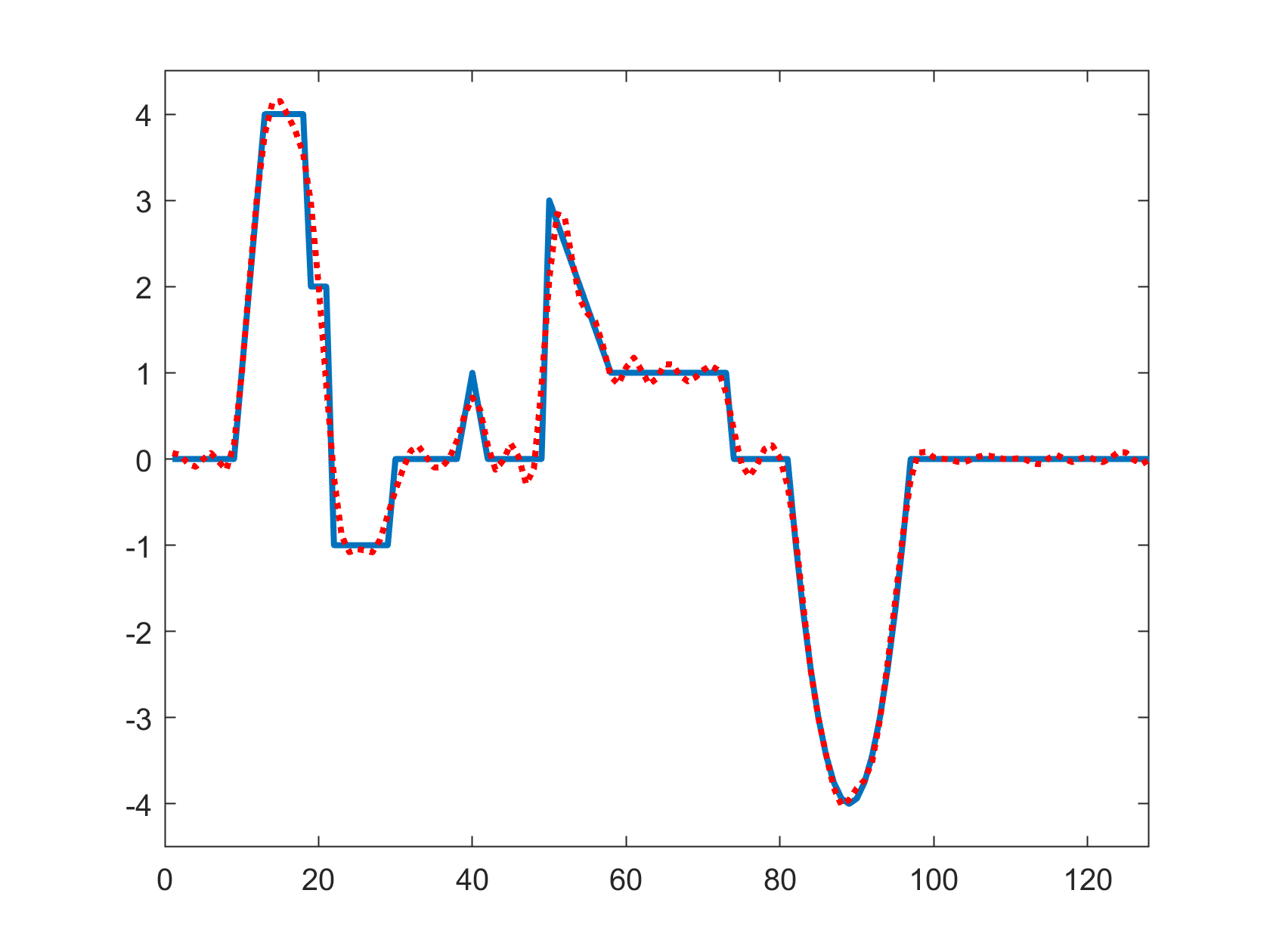}\\b)
		\end{minipage}
		\begin{minipage}{0.32\textwidth}
			\centering
			\includegraphics[width = \textwidth]{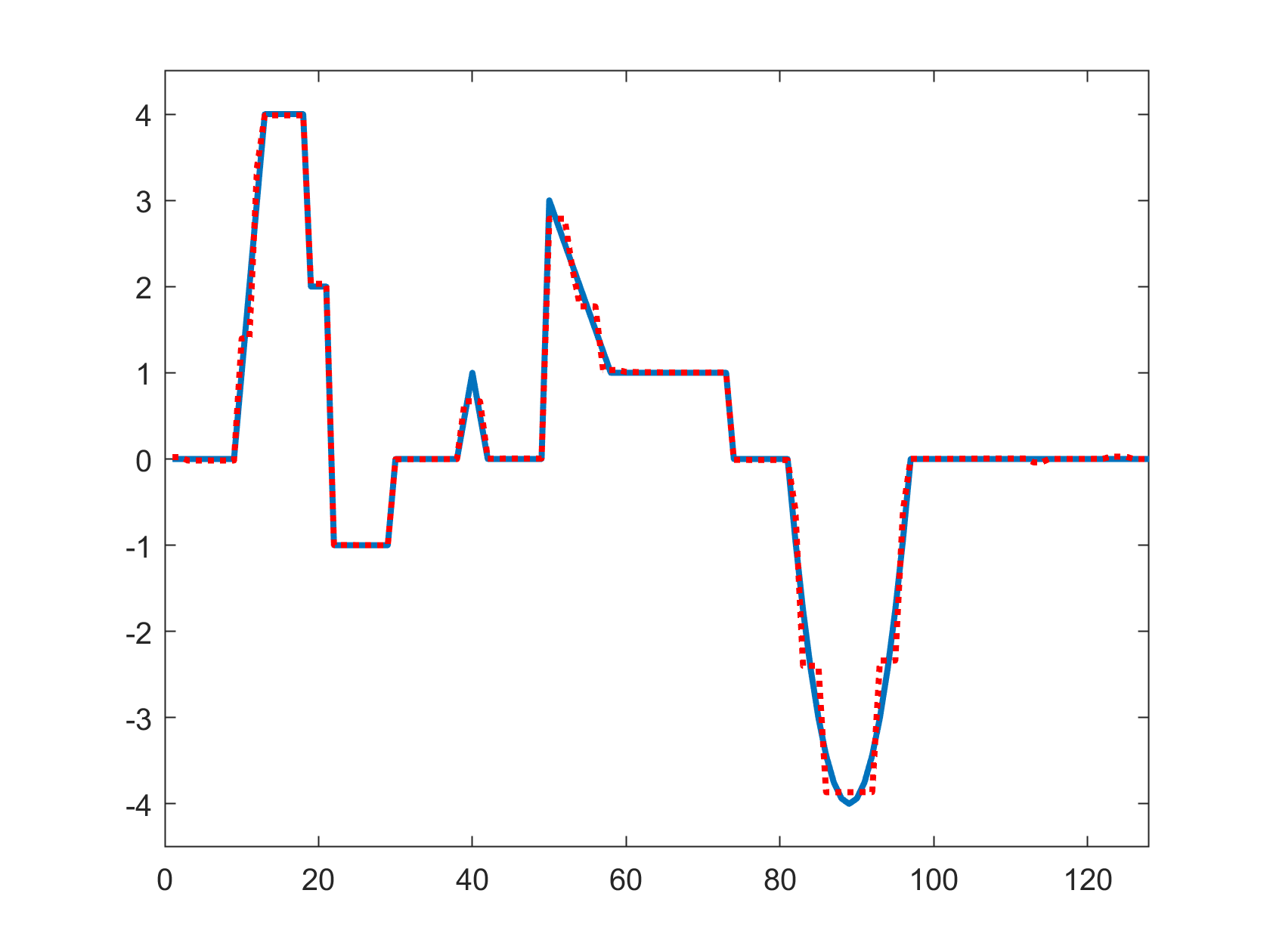}\\c)
		\end{minipage}
		\caption{Comparison of the use of a variety of $\ell_p$ norms to reconstruct a solution that has edges, using the exact value of the parameter $\by$, $\by_{\rm true}$. Blue line represents the true solution and the red line represents the blured and noisy data and the reconstructed solutions for different values of $p$. a) $\bx_{\rm true}$ and $\bd$  b) $p = 2$ and the RRE = 0.1396, c) $p = 0.5$ and RRE = 0.0974.}
		\label{fig: 1D Example}
	\end{figure}

	\subsection{One-dimensional deconvolution problem} Here we consider a one-dimensional convolution problem, where $\by$ is of dimension one. The forward model can be written as
	$$d(s)=\int_{a}^b g(s-s')x(s')ds',$$
	where we are are assuming a Gaussian kernel of the form 
	$$g(s)=\frac{1}{\sqrt{2\pi \sigma^2}} \rm{exp} \left( -\frac{s^2}{2\sigma^2}\right).$$
	To obtain the system $\bG\bx=\bd$ we consider $n$ discretization points and apply midpoint quadrature rule to approximate the integral. Furthermore, we consider zero boundary conditions to obtain a Toeplitz matrix. 
	We consider a small problem where the solution contains edges to demonstrate the need for the $p$-norm, and at the same time we can solve the $\ell_p$ problem easily without projecting onto a Krylov subspace, and compute the Jacobians by means of matrix decompositions.
	We use $n=128$, and the parameters $\sigma=2$. The condition number of $\bG(\by)$ is $1.7 \times 10^8$. The exact solution, represented by $\bx_{{\rm true}}$, is the vector of length $128$ shown in Figure \ref{fig: 1D Example}.The noise-free blurred signal, represented by $\bd_{{\rm true}}$, is computed as
	$\bd_{{\rm true}}=\bG\bx_{{\rm true}}$. The elements of the noise vector $\bepsilon$ are
	normally distributed with zero mean, and the standard deviation is
	chosen such that $\frac{\|\bepsilon\|_2}{\|\bd_{{\rm true}}\|_2}=0.01$. In this case, we say that the level of noise is $1\%$. The noisy right-hand side of our system is defined by $\bd=\bd_{{\rm true}}+\bepsilon$. Figure \ref{fig: 1D Example}) shows the solutions obtained when solving the regularized linear inverse problem with $\bG(\by_{{\rm true}})$ for $p=2$ and $p=0.05$. It is clear that for $p=0.05$ we can reconstruct some of the edges better than using $p=2$. The RRE values here can be used as a reference for the nonlinear problem as lower bounds since it would be hard to get better values than using $\by_{{\rm true}}$. For this example, we define 20 logarithmically spaced values of $\lambda$ from $10^{-4}$ to 1. We loop for all values of $\lambda$ and choose the $\lambda$ values that make RRE minimum for $\bx_{(i)}$.
	
\begin{figure}[h!]
		\centering
		\begin{minipage}{0.45\textwidth}
			\centering
			\includegraphics[width = \textwidth]{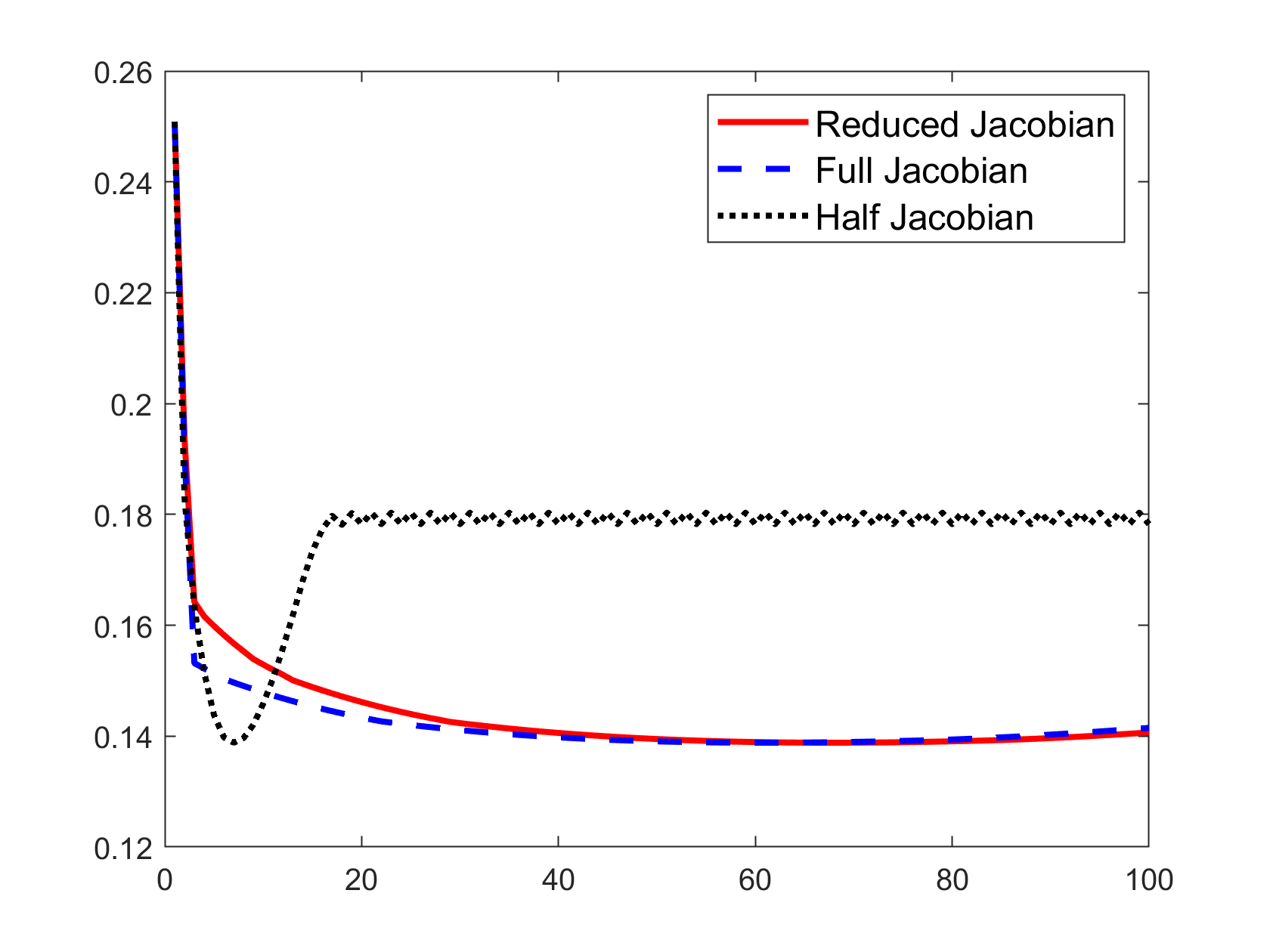}\\a)
		\end{minipage}
		\begin{minipage}{0.45\textwidth}
			\centering
			\includegraphics[width = \textwidth]{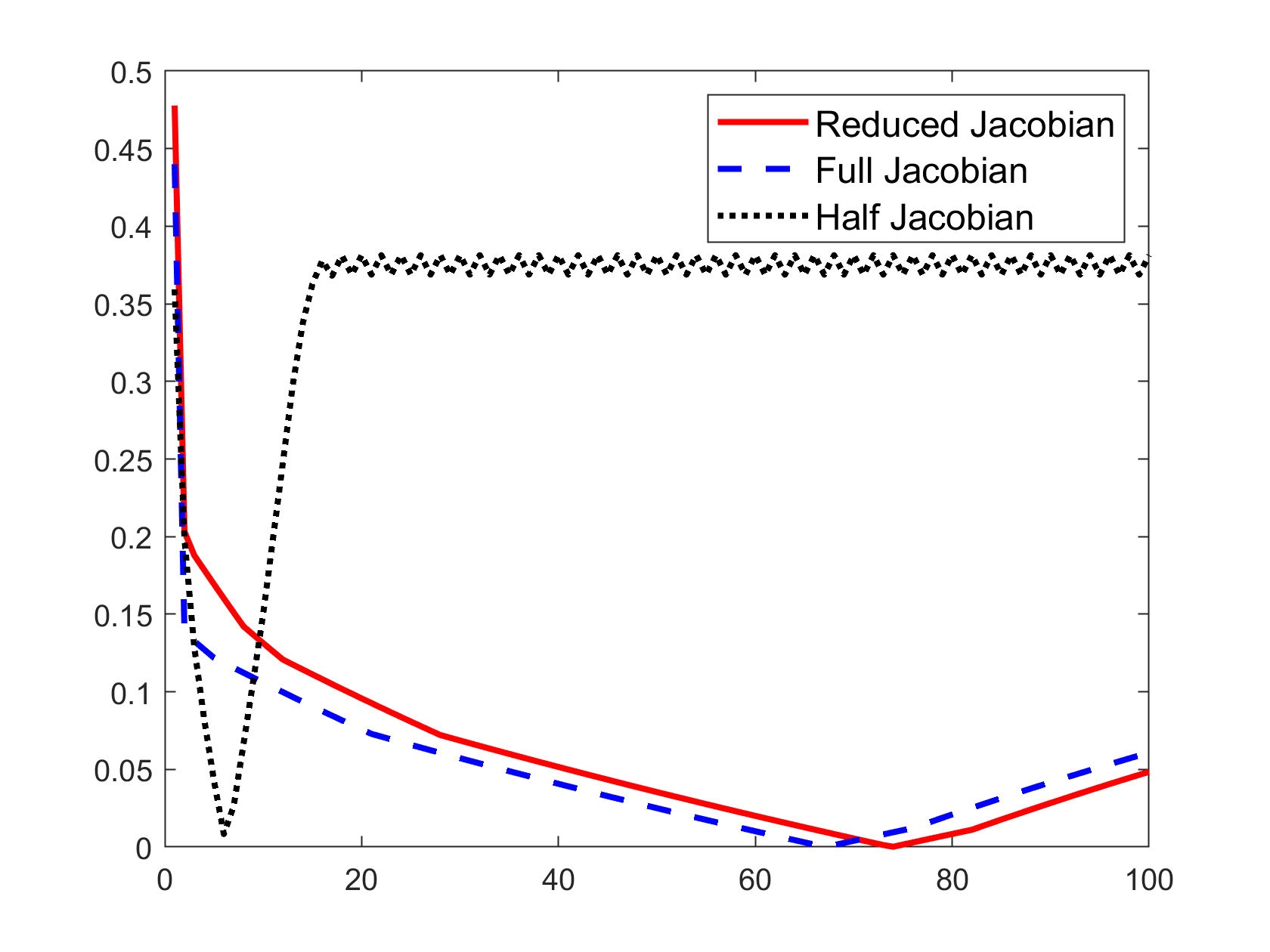}\\b)
		\end{minipage}
		\caption{1D example: Convergence curves to compare the used of different problems defined by the different Jacobians. a) RRE of the reconstructed image and b) RRE of the reconstructed parameters for 100 iterations of the VarPro method with reduced Jacobian (red solid line), full Jacobian (blue dashed line), and half Jacobian (black dotted line).}
		\label{fig: 1D-Convergence Curves}
	\end{figure}

	Figure \ref{fig: 1D-Convergence Curves} shows the convergence curves used to compare the different Jacobians for the case when $p=2$. 
	\begin{itemize}
		\item {\bf Reduced Jacobian}: We apply Algorithm \ref{Alg: Gen-VarPro} with steps 4-5 replaced by $\bs^{(i)} = \argmin_{\bs}\left\|\bJ^{(i)}\bs - \br^{(i)}\right\|^2_2$ with $\bJ^{(i)}=\frac{\partial \bG(\by)}{\partial \by_j} \bx^{(i)}$ and $\br^{(i)}=\bG(\by^{(i)})\bx^{(i)}-\bd$. This reduced Jacobian was used in \cite{chung2010efficient}.
		\item {\bf Full Jacobian}: We apply Algorithm \ref{Alg: Gen-VarPro} where the Jacobian contains $\bA_j$ and $\bB_j$ defined in \ref{eq: jthColumnJacobian}.
		\item {\bf Half Jacobian}: We apply Algorithm \ref{Alg: Gen-VarPro} where the Jacobian contains only the term $\bA_j$ defined in \ref{eq: jthColumnJacobian}. This was the approximation suggested in \cite{kaufman1975variable}.
	\end{itemize}
	
	We conclude that the reduced Jacobian and the full Jacobian makes the convergence behave similarly, but the reduced Jacobian is computational less expensive than the other two. They all produce similar RRE values (RRE($\bx^{(i)}$)=0.1388 and RRE($\by^{(i)}$)=0.008), but they achieved the minimum at different number of iterations (Reduced Jacobian: 68, Full Jacobian: 61, Half Jacobian: 6). All of them, have a semiconvergence behavior, so a stopping criteria is paramount. In the following large-scale numerical examples we will be testing $p\neq 2$ and we will use the reduced Jacobian.
	
	\subsection{Blind deconvolution model}
	We consider the problem of restoring images that  have been contaminated by blur and noise. The continuous space-invariant image deblurring problem can 
	be formulated as a Fredholm integral of the first kind, i.e., as an integral equation of the form 
	\begin{equation}\label{eq:im_debl}
		d(u,v)=\int_\Omega g(s-u,t-v)x(s,t)ds dt,\qquad (u,v)\in\Omega,
	\end{equation}
	where the function $d$ represents a blurred, but noise-free image, the function $x$ is the unknown image that we would like to recover, and $g$ is a smooth kernel with compact support. The integral 
	operator in \eqref{eq:im_debl} is compact. Therefore, the solution of \eqref{eq:im_debl} 
	is an ill-posed problem. Discretizing \eqref{eq:im_debl} gives the system
	$\bG \bx = \bd$, with the matrix $\bG$ being the sum of a block Toeplitz with Toeplitz block
	matrix and a correction of small norm due to the boundary conditions that are imposed; 
	see, e.g., \cite{hansen2006deblurring} for more details on image deblurring. 
	
	As a first step towards the numerical examples, we describe the blind deconvolution problem that is obtained by the convolution equation
	\begin{equation}\label{eq: conv}
		\tilde d(s, t) = p(s,t) \star x(s,t) + \epsilon(s, t),
	\end{equation}
	where the true image $x(s, t)$ is convolved with the point spread function (PSF) $p(s,t)$ and noised with the additive noise $\epsilon(s,t)$ to obtain the blurred and noisy image $\tilde d(s,t)$. The blur is assumed to be spatially invariant, that is shift invariant. The case when the operator that models the forward problem is not known exactly is referred to as blind deconvolution.  We consider the discretized version of the model of the image formation where 
	\begin{itemize}
		\item $\bx$ is the vector that represent the unknown true image that we aim to reconstruct.
		\item $\bd$ is the vector that represents the measurements that are blurred and noised.
		\item $\bG(\by)$ is an ill-conditioned matrix that describes the parameter-to-observable map representing the PSF that is defined by the parameter's vector $\by$. Usually $\bG$ is a sparse and/or structured matrix.
		\item $\by$ is the vector of parameters that define the blurring operator. We assume that the vector $\by$ contains a small set of parameters that define the PSF. In our case, it will represent a Gaussian PSF where we have
		\begin{equation}\label{eq: PSFfo}
			p(s,t) = \frac{1}{2\pi \sqrt{\delta}}\exp\left( -\frac{1}{2} 
			\begin{bmatrix} 
				s\\
				t
			\end{bmatrix}^T 
			\begin{bmatrix}
				\sigma_1^2 & \rho^2 \\ \rho^2  & \sigma_2^2
			\end{bmatrix}^{-1}
			\begin{bmatrix}
				s\\
				t
			\end{bmatrix}
			\right),
		\end{equation}
		where $\delta = \sigma_1^2\sigma_2^2 - \rho^4 > 0$. This recommends that the set of parameters can be represented by three values $\sigma_1$, $\sigma_2$, and $\rho$.
		\item $\bepsilon$ represents the noise present in the measurements $\bd$ that may stem from rounding errors or measurement inaccuracies. 
	\end{itemize}
	
	In this section we use a special image deblurring problem that can be modeled as a separable nonlinear inverse problem to illustrate the performance of the approach that we propose. The matrix that models the blurring operator is defined by $\bG(\by) = \bG (\bP(\by))$, where $\bP(\by)$ is the PSF (see Figure \ref{fig: Satellite_problemSetup} for an example). We scale the entries of the PSF to sum up to 1.
	The reduced Jacobian can be computed using the chain rule. Specifically, 
	\begin{equation}\label{eq: Jacobian}
		\bJ = \frac{\partial}{\partial \by}[\bG (\bP(\by))\bx] \overset{(a)}{=} \frac{\partial}{\partial\bP}[\bG(\bP(\by))\bx]\frac{\partial}{\partial \by}{\bP(\by)} \overset{(b)}{=} \bG(\bX)\frac{\partial}{\partial \by}{\bP(\by)} \overset{(c)}{=} \bG(\frac{\partial}{\partial \by}{\bP(\by)})\bx.
	\end{equation}
	The (a) equality in \eqref{eq: Jacobian} follows from the chain rule, (b) follows from the commutative property of the convolution operator, that is $\bG(\bP(\by))\bx = \bG(\bX)\rm{vec}(\bP(\by))$, and c) follows by the same argument as in (b). As mentioned in \cite{chung2010efficient}, the Jacobians can also be computed using finite differences. 
	
	The task of determining both a restored image and an improved approximation of the blurring matrix is commonly referred to as blind deconvolution. There is quite a bit of work done in the blind deconvolution problem (see \cite{chan1998total, chan2000convergence, desidera2006application, mugnier2004mistral, yu2012blind} for some selected references).
	\subsection{Regularization models that we consider}
	For our numerical experiments we consider three main cases for the regularization term $\|\bL\bx\|_p^p$ in \eqref{eq: lp}. 
	\paragraph{Case 1} As a first model we use $\bL = \bI$ because we would like to reconstruct a solution $\bx$ whose elements have minimal $p$-norm. 
	\paragraph{Case 2} As a second model we use a two-level framelet analysis operator $\bL$ since it is well-known that images have sparse representation in the framelet domain. We recall that framelets are extensions of 
	wavelets. We define the framelets that we use as follows.
	\begin{definition}\label{def: Tightframe}
		Let $\bW\in\R^{r\times n}$ with $1\leq n\leq r$. The set of the rows of $\bW$ is a framelet system
		for $\R^n$ if for all $\bx\in\R^n$ it holds that
		\begin{equation} \label{eq: tightframe}
			\|\bx\|_2^2=\sum_{j=1}^{r}{(\bw\t_j\bx)^2},
		\end{equation}
		where $\bw_j\in\R^n$ denotes the $j$th row of the matrix $\bW$ (written as a column vector), 
		i.e., $\bW=[\bw_1,\bw_2,\ldots,\bw_r]\t$. The matrix $\bW$ is referred to as an analysis operator 
		and $\bW\t$ as a synthesis operator. 
	\end{definition}
	We use the same tight frames as in \cite{buccini2020modulus, buccini2, COS09b}, i.e., the system of linear 
	B-splines. This system is formed by a low-pass filter $\bW_{0}\in\mathbb{R}^{n \times n}$ 
	and two high-pass filters $\bW_1,\bW_2\in\mathbb{R}^{n \times n}$, whose corresponding masks 
	are
	\begin{equation*}
		\bw^{(0)}=\frac{1}{4}[1,2,1], \quad \bw^{(1)}=\frac{\sqrt{2}}{4}[1,0,-1], 
		\quad \bw^{(2)}=\frac{1}{4}[-1,2,-1].
	\end{equation*}
	The analysis operator $\bW$ in one space-dimension is derived from these masks and by 
	imposing reflexive boundary conditions to ensure that $\bW\t\bW=\bI$. The so-determined
	filter matrices are 
	\[ 
	\bW_{0}=\frac{1}{4}\begin{bmatrix}
		3 & 1 & 0& \dots &0\\
		1 & 2 & 1  \\
		&  \ddots& \ddots & \ddots \\
		& &  1 & 2 & 1 \\
		0& \dots & 0 & 1 &3
	\end{bmatrix},\quad
	\bW_{1}=\frac{\sqrt2}{4}\begin{bmatrix}
		-1 & 1 & 0& \dots &0\\
		-1 & 0 & 1  \\
		&  \ddots& \ddots & \ddots \\
		& &  -1 & 0 & 1 \\
		0& \dots & 0 & -1 &1
	\end{bmatrix},
	\]
	and 
	\[
	\bW_{2}=\frac{1}{4}\begin{bmatrix}
		1 & -1 & 0& \dots &0\\
		-1 & 2 & -1  \\
		&  \ddots& \ddots & \ddots \\
		& &  -1 & 2 & -1 \\
		0& \dots & 0 & -1 &1
	\end{bmatrix}.
	\]
	The corresponding two-dimensional operator $\bW$ is given by 
	\begin{equation}\label{eq: Tightframe_W}
		\bW=\begin{bmatrix}
			\bW_{0}\otimes \bW_{0}\\
			\bW_{0}\otimes \bW_{1}\\
			\bW_{0}\otimes \bW_{2}\\
			\bW_{1}\otimes \bW_{0}\\
			\vdots\\
			\bW_{2}\otimes \bW_{2}
		\end{bmatrix},
	\end{equation}
	where $\otimes$ denotes the Kronecker product. This matrix is not explicitly formed. We
	note that the evaluation of matrix-vector products with $\bW$ and $\bW\t$ is inexpensive, 
	because the matrix $\bW$ is very sparse.
	\paragraph{Case 3} As a third model we use a regularization matrix that represents discretizations of the first and second derivative operator. Let 
	\begin{equation}\label{D1}
		\bL_{1,1}= 
		\begin{bmatrix} 
			1 &-1 & & & &0 \\
			&1 &-1 & & & \\
			& &.  &. &. & \\
			0& &  & &1  &-1 
		\end{bmatrix}\in \R^{(n-1)\times n},
	\end{equation}
	be a scaled discretization of the first derivative operator at equidistant points 
	in one space-dimension with periodic boundary conditions.
	$\bL_1 = \bL_{1,1}\otimes \bI + \bI \otimes \bL_{1,1}$.
	Similarly, we define $\bL_2 = \bL_{2,2}\otimes \bI + \bI \otimes \bL_{2,2}$, where 
	\begin{equation}\label{Lex}
		\bL_{2,2}= 
		\begin{bmatrix} 
			-1 &2 &-1 & & & \\
			&. &.  &. & & \\
			& &.  &. &. & \\
			& &  &-1 &2  &-1 
		\end{bmatrix}\in\R^{(n-2)\times n},
	\end{equation}
	
	\subsection{Example 1: Satellite test problem}
\begin{figure}[h!]
		\centering
		\begin{minipage}{0.3\textwidth}
			\centering
			\includegraphics[width = \textwidth]{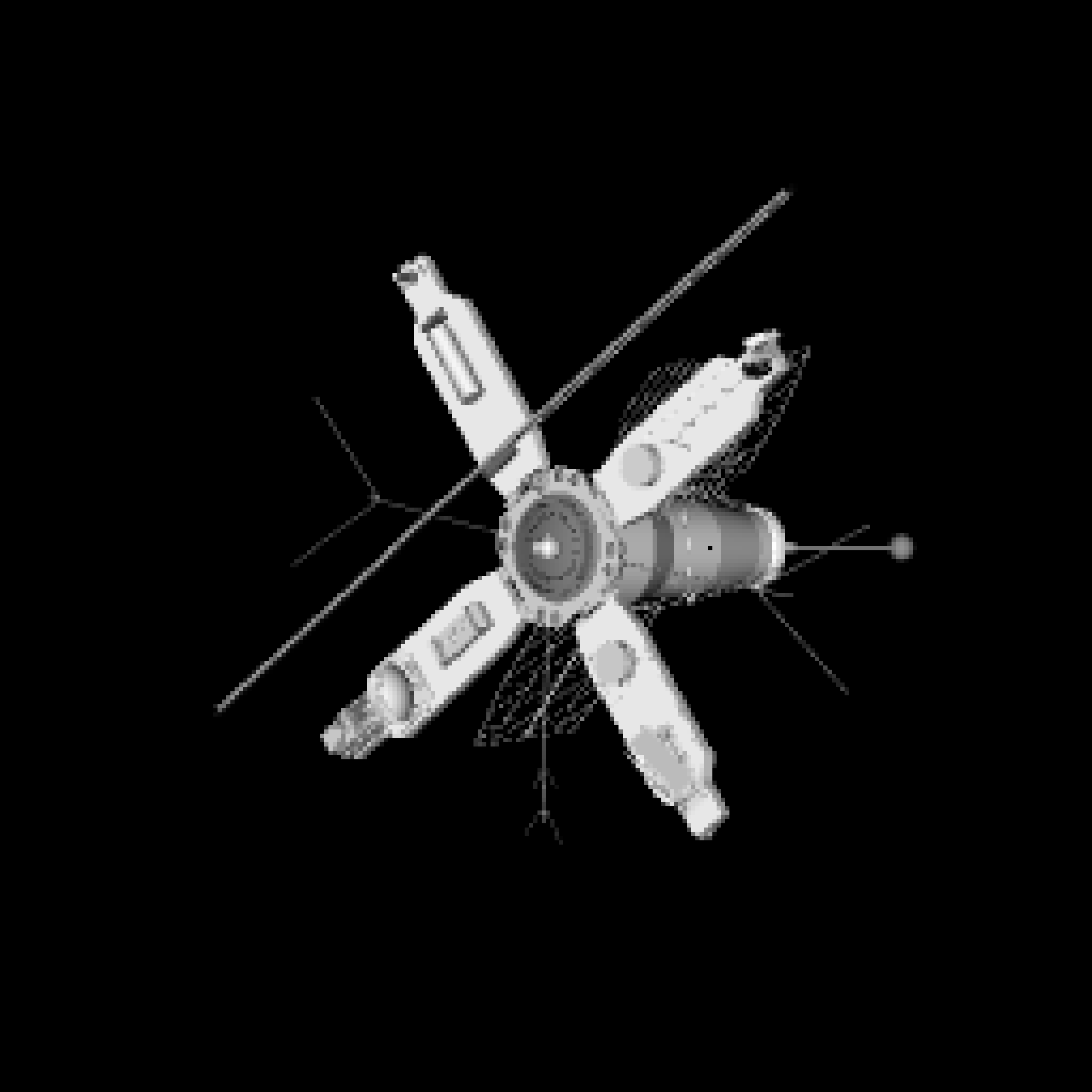}\\(a)
		\end{minipage}
		\begin{minipage}{0.3\textwidth}
			\centering
			\includegraphics[width = \textwidth, height = \textwidth]{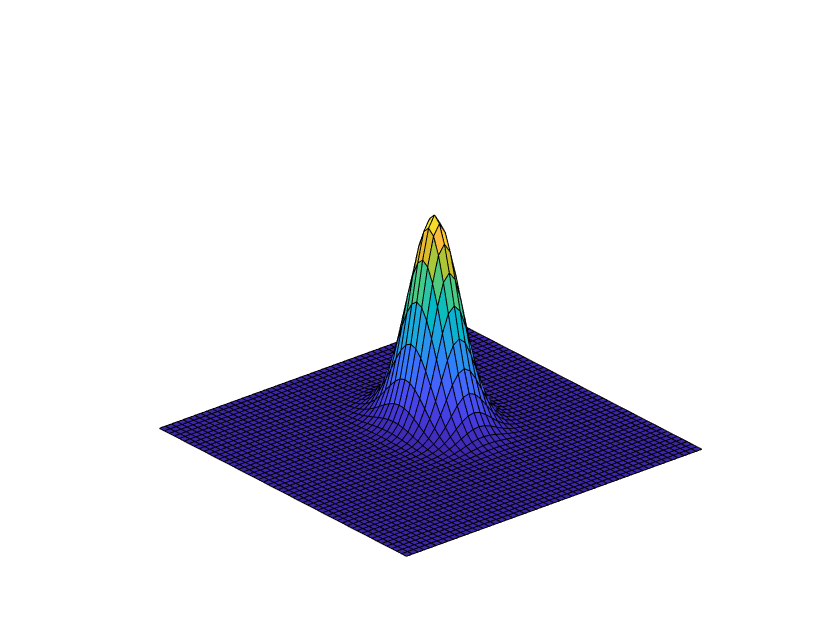}\\(b)
		\end{minipage}
		\begin{minipage}{0.3\textwidth}
			\centering
			\includegraphics[width = \textwidth]{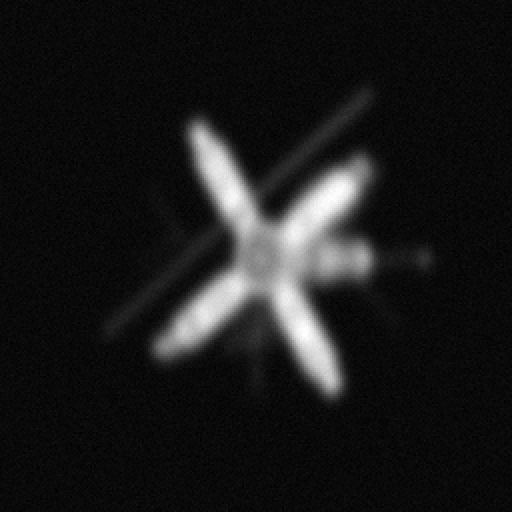}\\(c)
		\end{minipage}
		\caption{Satellite test problem setup. a) The true image, b) The PSF, c) Blurred and noisy image with $5\%$ Gaussian white noise.}
		\label{fig: Satellite_problemSetup}
	\end{figure}

	The first example considers the satellite test problem from \cite{gazzola2018ir}. The exact image is shown in Figure \ref{fig: Satellite_problemSetup}(a) and has size $256\times 256$. We blur the image with a Gaussian PSF shown in Figure \ref{fig: Satellite_problemSetup}(b) that is defined by the true parameters $\by_{\rm true} = (\sigma_1, \sigma_2, \rho) = (1.5,2.0,1.0)$ and the initial approximation of the parameters to be $\by = (3.0, 4.0, 2.0)$. The blurred image with $5\%$ Gaussian noise is shown in Figure \ref{fig: Satellite_problemSetup}(c). 
	\subsubsection{Investigation of the quality of the reconstructed image $\bx$ and parameters $\by$} 
\begin{figure}[h!]
		\centering
		\begin{minipage}{0.24\textwidth}
			\centering
			\includegraphics[width = \textwidth]{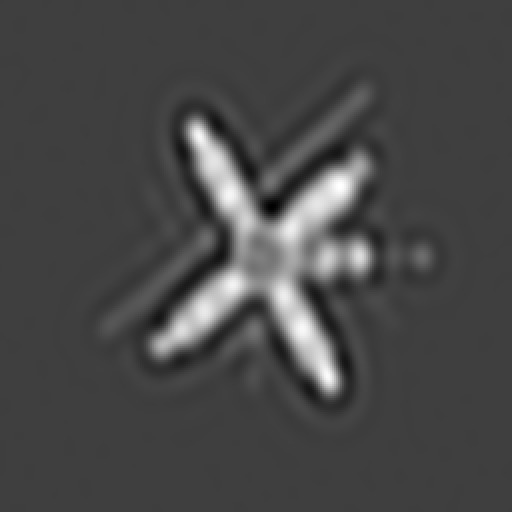}\\(a)
		\end{minipage}
		\begin{minipage}{0.24\textwidth}
			\centering
			\includegraphics[width = \textwidth]{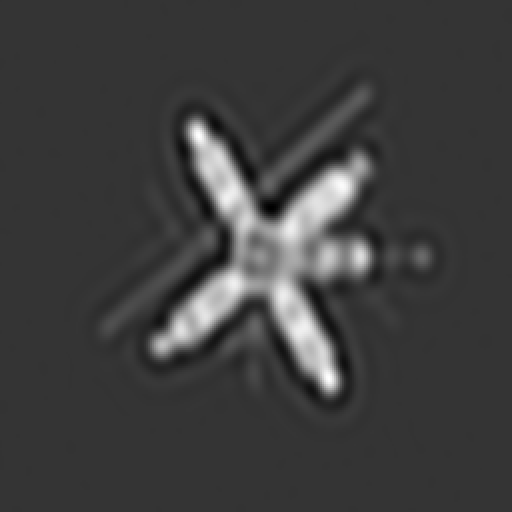}\\(b)
		\end{minipage}
		\begin{minipage}{0.24\textwidth}
			\centering
			\includegraphics[width=\textwidth]{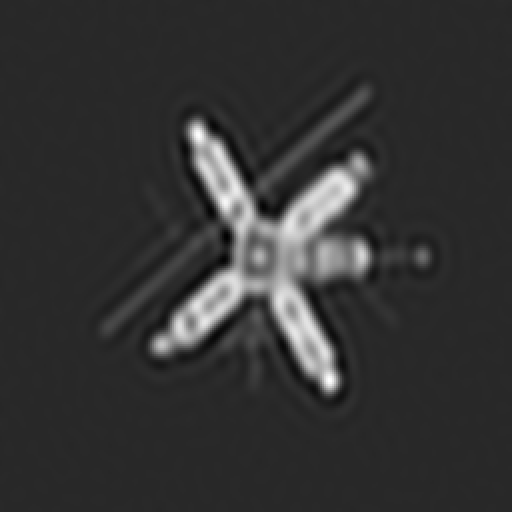}\\(c)
		\end{minipage}
		\begin{minipage}{0.24\textwidth}
			\centering
			\includegraphics[width = \textwidth]{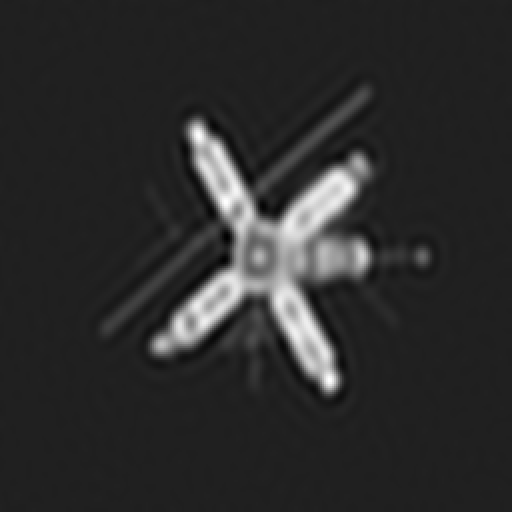}\\(d)
		\end{minipage}\\
		\begin{minipage}{0.24\textwidth}
			\centering
			\includegraphics[width = \textwidth]{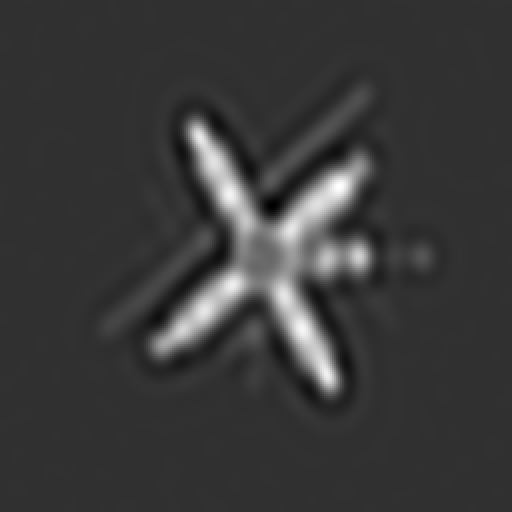}\\(e)
		\end{minipage}
		\begin{minipage}{0.24\textwidth}
			\centering
			\includegraphics[width = \textwidth]{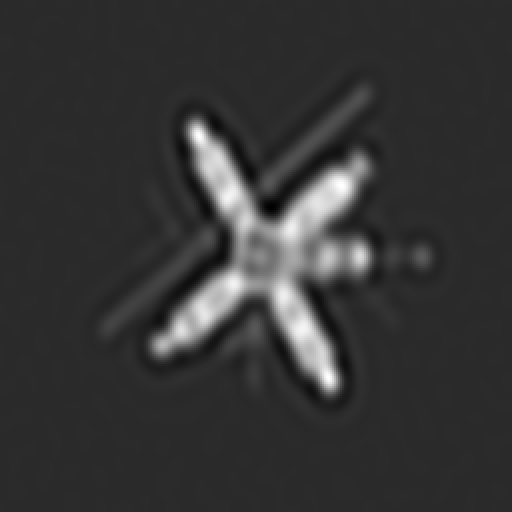}\\(f)
		\end{minipage}
		\begin{minipage}{0.24\textwidth}
			\centering
			\includegraphics[width=\textwidth]{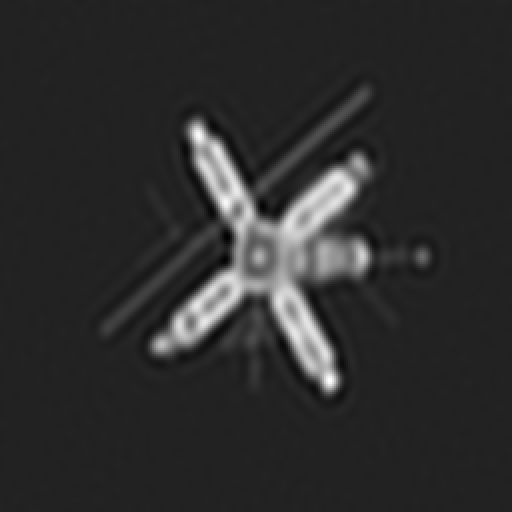}\\(g)
		\end{minipage}
		\begin{minipage}{0.24\textwidth}
			\centering
			\includegraphics[width = \textwidth]{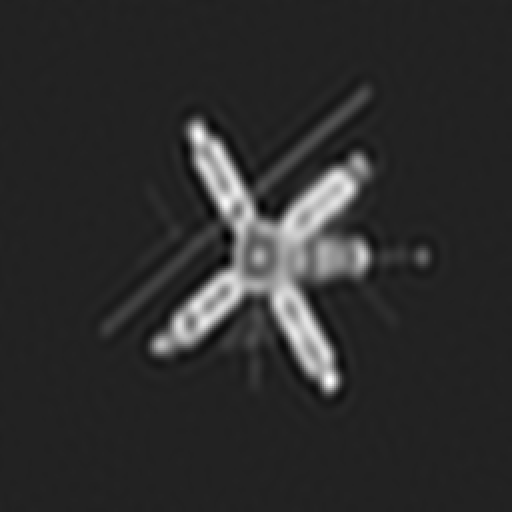}\\(h)
		\end{minipage}
		\caption{Reconstructions obtained by NLS \cite{chung2010efficient} and NLS$_{\ell_{p}-\ell_{q}}$. Panels (a-d) show reconstructions using NLS \cite{chung2010efficient} and panels (e-h) show reconstructions using NLS$_{\ell_{p}-\ell_{q}}$ at iterations $i = 1, 3, 7, 11$ respectively. }
		\label{Figure: satellite_rec}
	\end{figure}

\begin{figure}[h!]
		\centering
		\begin{minipage}{0.3\textwidth}
			\centering
			\includegraphics[width = \textwidth]{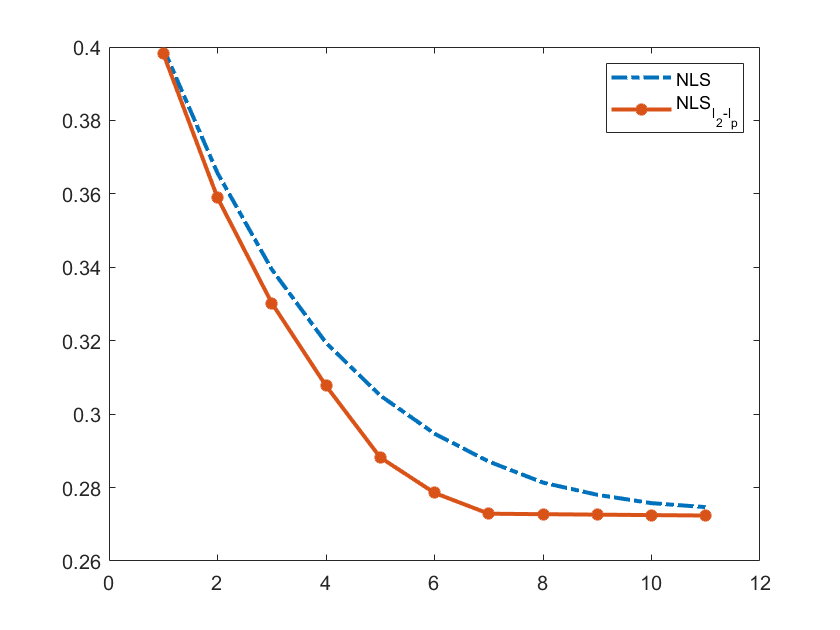}\\(a)
		\end{minipage}
		\begin{minipage}{0.3\textwidth}
			\centering
			\includegraphics[width = \textwidth]{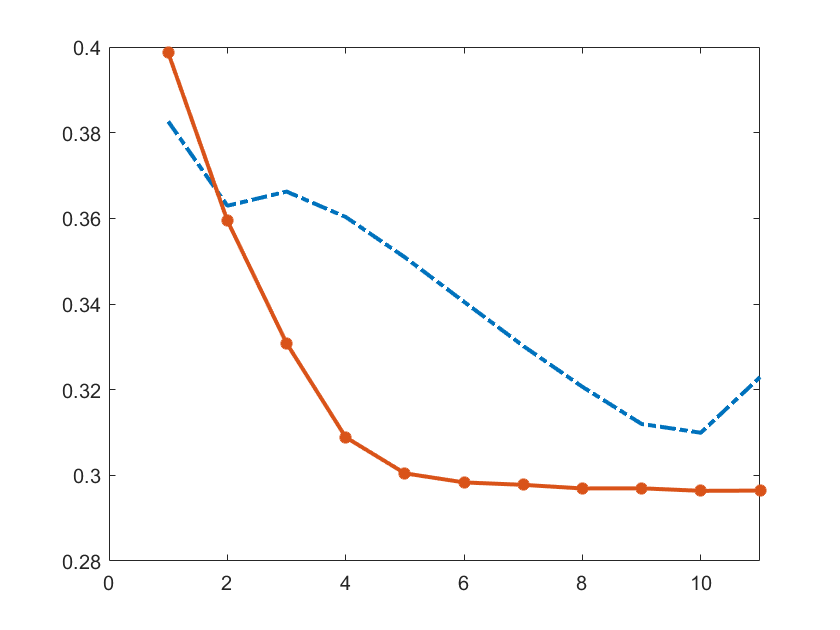}\\(b)
		\end{minipage}
		\begin{minipage}{0.3\textwidth}
			\centering
			\includegraphics[width = \textwidth]{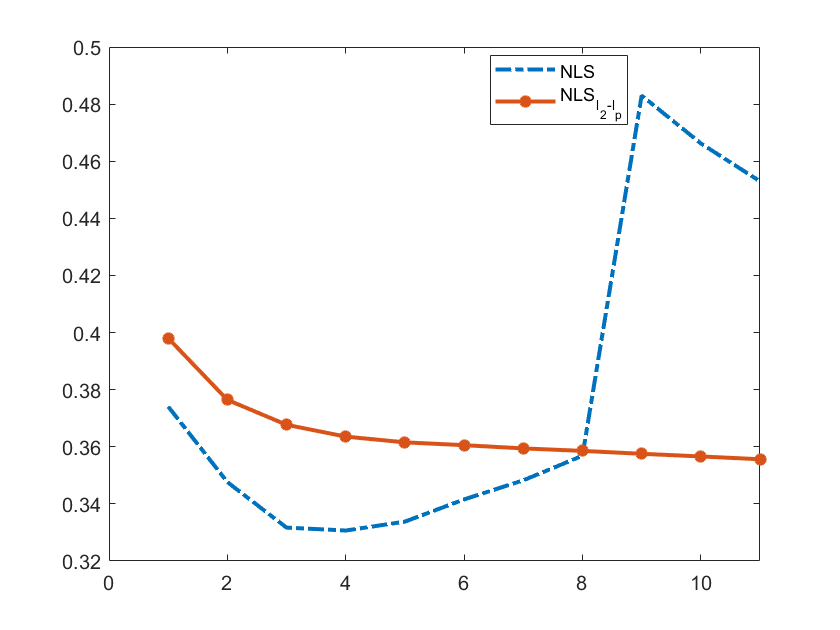}\\(c)
		\end{minipage}
		\caption{RRE versus the number of iterations. The dotted red line represents RRE for NLS and stared blue line represents the RRE for  NLS$_{\ell_{2}-\ell_{q}}$ for different noise levels added in the available data. a) 1\% Gaussian noise, b) 5\% Gaussian noise, c) 10\% Gaussian noise.}
		\label{fig: satellite}
	\end{figure}

In Figure \ref{Figure: satellite_rec}, the reconstruction obtained using $p = 2$ (NLS) and $p = 1.1$ (NLS$_{\ell_{2}-\ell_{p}}$) are presented. We chose $p = 1.1$ to make sure we are solving a strictly convex optimization for $\bx$ (step 3 in Algorithm \ref{Alg: lp VarPro}). Nevertheless, it is known in image deblurring that when the desired solution is sparse, choosing $p<1$ enhances the quality of the reconstructed solution yet the problem to be solved may be non-convex. For a fair comparison, we let $\bL = \bI$, however, we acknowledge that the results presented here can be improved by different choices of $\bL \neq \bI$. We vary the noise level from $1\%$ to $5\%$ to illustrate the robustness of the method with respect to noise. Indeed, in the RRE of the desired parameters $\by$ shown in Figure \ref{fig: satellite}, we observe a quicker decay, as well as the values that stabilize after a certain number of iterations. A similar behaviour is observed for large noise levels, as shown in Figure \ref{fig: satellite}(c). In addition, we investigate the objective function values and the RRE of the reconstructed image. While we observe slightly improvements in the values of the objective function when NLS$_{\ell_2-\ell_p}$ is used, we obtain faster convergence and smaller RREs (see Figure \ref{fig: satellite_objfunction}(a) and (b)).
	More insight on the convergence of the relative function value (RelFuncValue), relative gradient norm (RelGradNorm), RRE of the parameters, and the RRE of the image are shown in Tables \ref{Tab: satellite_Hybrid_info} and \ref{Tab: satellite_lp} for iterations $i =1,2,\dots,11$.
	
	\subsubsection{Investigation on different choices of $p$} 
\begin{figure}[h!]
		\centering
		\begin{minipage}{0.45\textwidth}
			\centering
			\includegraphics[width = \textwidth]{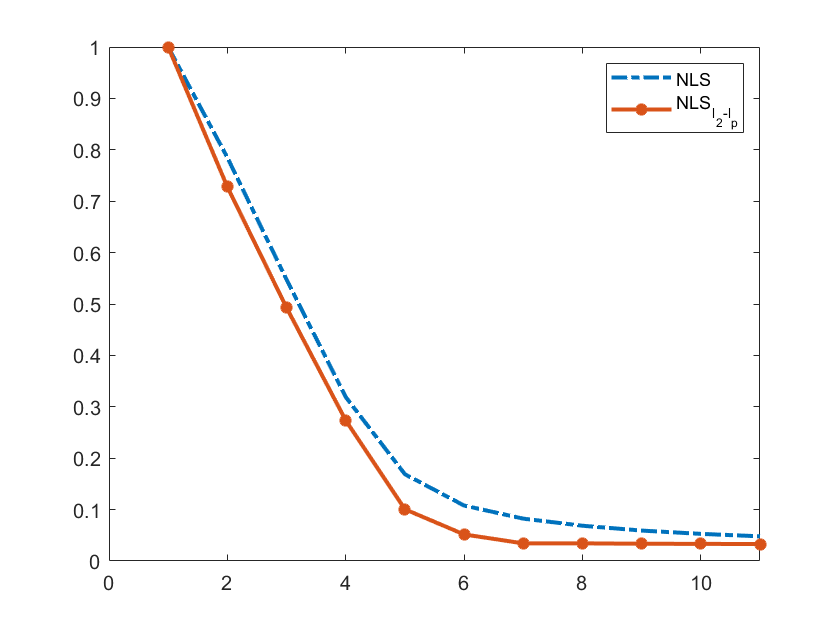}\\(a)
		\end{minipage}
		\begin{minipage}{0.45\textwidth}
			\centering
			\includegraphics[width =\textwidth]{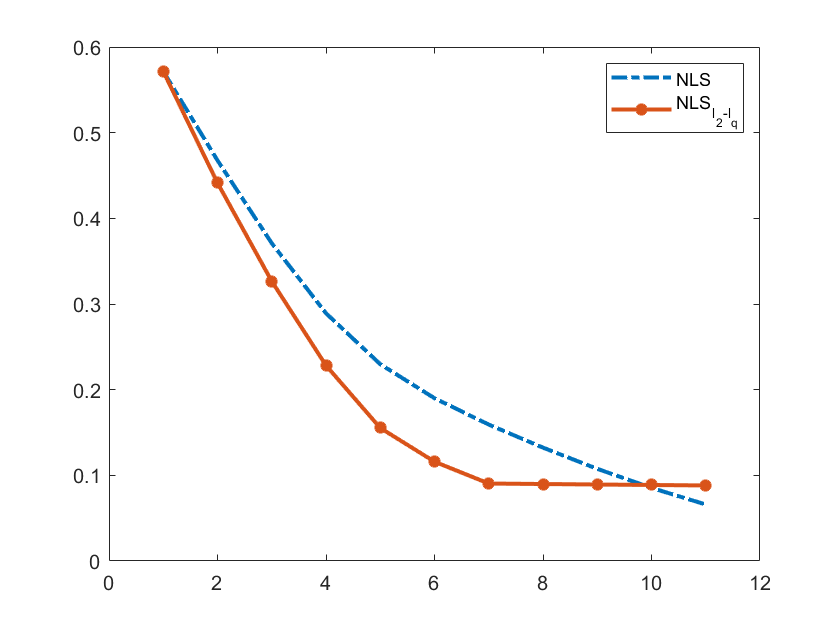}\\(b)
		\end{minipage}
		\caption{Objective function behaviour. a) RRE of the reconstructed image vs. the number of iterations. b) RRE of the reconstructed parameters vs. the number of iterations.}
		\label{fig: satellite_objfunction}
	\end{figure}

In this subsection we consider the true object to be the satellite test image, blurred with Gaussian PSF defined by the parameters $\by_{\rm true} = (\sigma_1, \sigma_2, \rho) = (1.5, 2.0, 0.5)$ and perturbe the available measurements with $1\%$ Gaussian white noise. We vary the parameter $p$ to investigate the behaviour of the quality of the reconstructed image $\bx$ and the reconstructed parameters $\by$. We choose the starting approximation for the parameters to be $\by^{(0)} = (3.0, 4.0, 1.5)$. 
	Here we choose the regularization matrix $\bL$ to be the discretization of the first derivative operator $\bL = \bL_1$ that enforces extra properties on the desired solution $\bx$, hence we seek to reconstruct solutions such that their derivatives are sparse when $0<p<2$. For our investigation, we set $p$ to be 0.8, 1.2, and 2.0. 
	We present the objective function behaviour and the convergence of the RRE of the parameters $\by$ in figure \eqref{fig: satellite_objfunction}(a) and (b) respectively.
	
	\begin{table}[h!]
		\centering
		\caption{Convergence of NLS \cite{chung2010efficient} method for Satellite test problem with 1\% white Gaussian noise.}
		\label{Tab: satellite_Hybrid_info}
		\begin{tabular}{l|llll}\hline
			i & RelFuncValue & RelGradNorm & RRE($\by^{(i)}$) & RRE($\bx^{(i)}$) \\\hline
			1.0     & 1.0000       & 1.0000      & 0.5716       & 0.3998      \\
			2.0     & 0.7850       & 0.8753      & 0.4676       & 0.3657      \\
			3.0     & 0.5476       & 0.7190      & 0.3709       & 0.3394      \\
			4.0     & 0.3190       & 0.5313      & 0.2892       & 0.3194      \\
			5.0     & 0.1688       & 0.3636      & 0.2297       & 0.3051      \\
			6.0     & 0.1079       & 0.2761      & 0.1900       & 0.2947      \\
			7.0     & 0.0824       & 0.2333      & 0.1592       & 0.2871      \\
			8.0     & 0.0686       & 0.2087      & 0.1323       & 0.2814      \\
			9.0     & 0.0592       & 0.1874      & 0.1076       & 0.2780      \\
			10.0    & 0.0528       & 0.1722      & 0.0853       & 0.2758      \\
			11.0    & 0.0480       & 0.1601      & 0.0660       & 0.2747   \\\hline  
		\end{tabular}
	\end{table}
	
	\begin{table}[h!]
		\centering
		\caption{Convergence of NLS$_{\ell_{2}-\ell_{q}}$ method for Satellite test problem with 1\% white Gaussian noise.}
		\label{Tab: satellite_lp}
		\begin{tabular}{l|llll}\hline
			i & RelFuncValue & RelGradNorm & RRE($\by^{(i)}$) & RRE($\bx^{(i)}$) \\\hline
			1.0     & 1.0000       & 1.0000      & 0.5716       & 0.3982      \\
			2.0     & 0.7278       & 0.8238      & 0.4413       & 0.3590      \\
			3.0     & 0.4929       & 0.6718      & 0.3266       & 0.3303      \\
			4.0     & 0.2735       & 0.4981      & 0.2282       & 0.3077      \\
			5.0     & 0.1006       & 0.2817      & 0.1553       & 0.2882      \\
			6.0     & 0.0517       & 0.1791      & 0.1159       & 0.2785      \\
			7.0     & 0.0343       & 0.1244      & 0.0904       & 0.2729      \\
			8.0     & 0.0340       & 0.1231      & 0.0899       & 0.2727      \\
			9.0     & 0.0336       & 0.1218      & 0.0893       & 0.2726      \\
			10.0    & 0.0333       & 0.1204      & 0.0887       & 0.2725      \\
			11.0    & 0.0330       & 0.1193      & 0.0882       & 0.2724     \\\hline
		\end{tabular}
	\end{table}
	
	\begin{figure}[h!]
		\centering
		\begin{minipage}{0.45\textwidth}
			\centering
			\includegraphics[width = \textwidth]{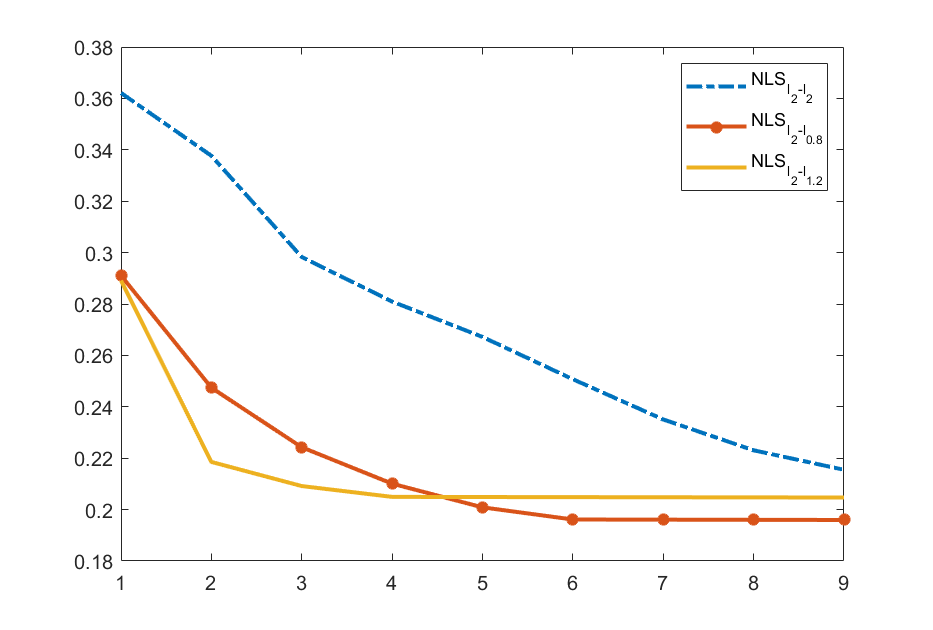}\\(a)
		\end{minipage}
		\begin{minipage}{0.45\textwidth}
			\centering
			\includegraphics[width =\textwidth]{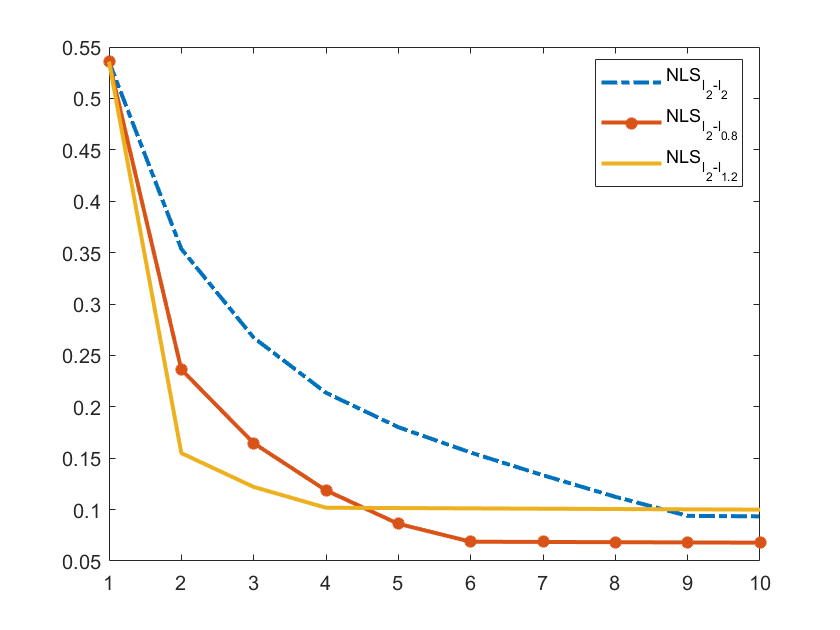}\\(b)
		\end{minipage}
		\caption{Satellite test problem-investigation of choices of $p$. a) RRE of the reconstructed image $\bx$ in the vertical axis vs. the number of iterations in the horizontal axis. b) RRE of the reconstructed parameters $\by$ vs. the number of iterations.}
		\label{fig: satellite_objfunction_2}
	\end{figure}
	
	\subsection{Example 2: Grain test problem} 
\begin{figure}[h!]
		\centering
		\begin{minipage}{0.3\textwidth}
			\centering
			\includegraphics[width = \textwidth]{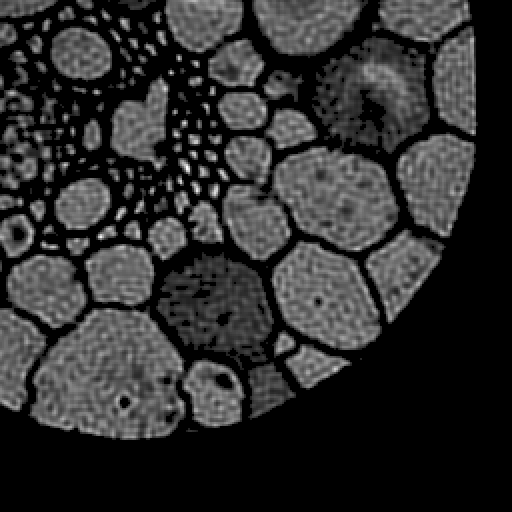}\\(a)
		\end{minipage}
		\begin{minipage}{0.3\textwidth}
			\centering
			\includegraphics[width = \textwidth, height = \textwidth]{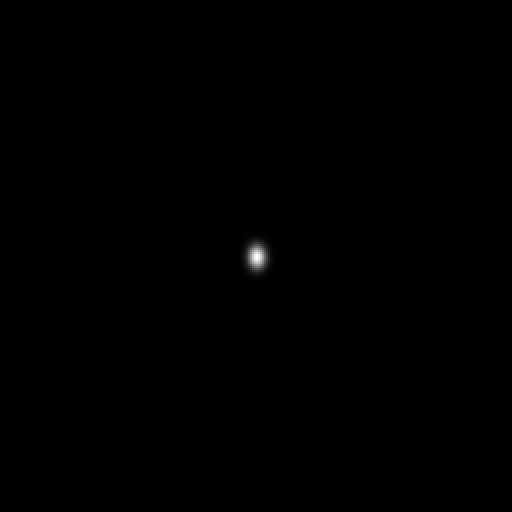}\\(b)
		\end{minipage}
		\begin{minipage}{0.3\textwidth}
			\centering
			\includegraphics[width = \textwidth]{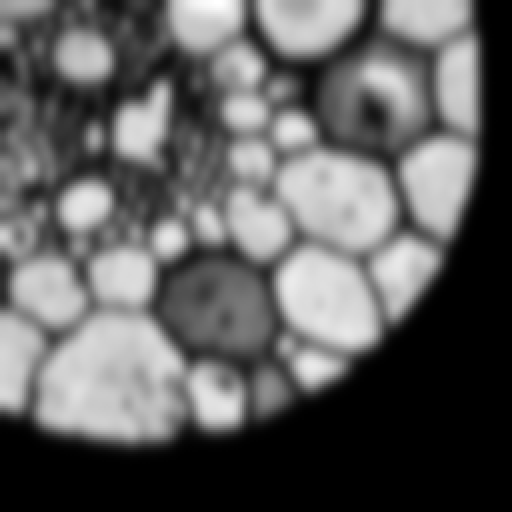}\\(c)
		\end{minipage}
		\caption{Grain test problem setup. a) The true image, b) the PSF, and c) Blurred and noisy image with $1\%$ Gaussian white noise.}
		\label{fig:grain_problemSetup}
	\end{figure}

\begin{figure}[h!]
		\centering
		\begin{minipage}{0.24\textwidth}
			\centering
			\includegraphics[width = \textwidth]{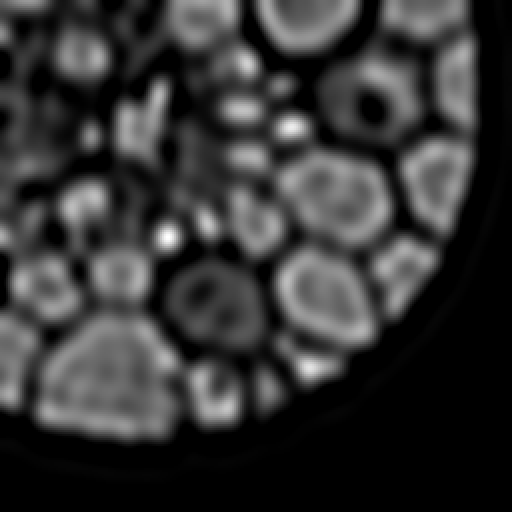}\\(a)
		\end{minipage}
		\begin{minipage}{0.24\textwidth}
			\centering
			\includegraphics[width = \textwidth]{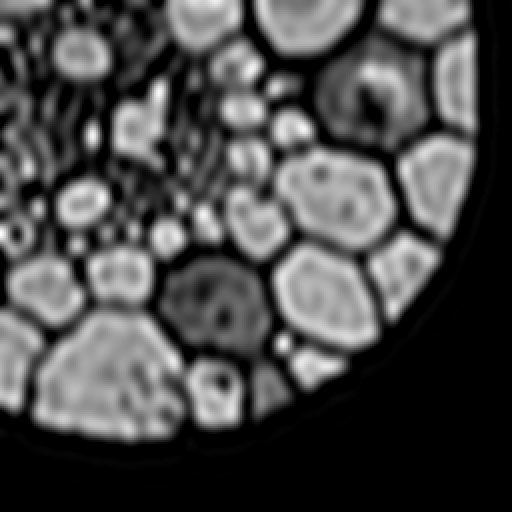}\\(b)
		\end{minipage}
		\begin{minipage}{0.24\textwidth}
			\centering
			\includegraphics[width=\textwidth]{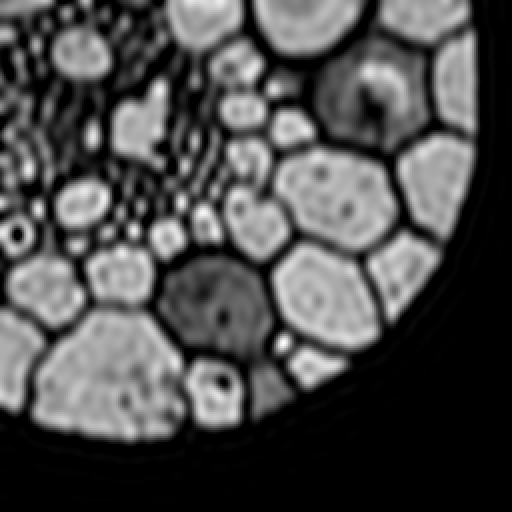}\\(c)
		\end{minipage}
		\begin{minipage}{0.24\textwidth}
			\centering
			\includegraphics[width = \textwidth]{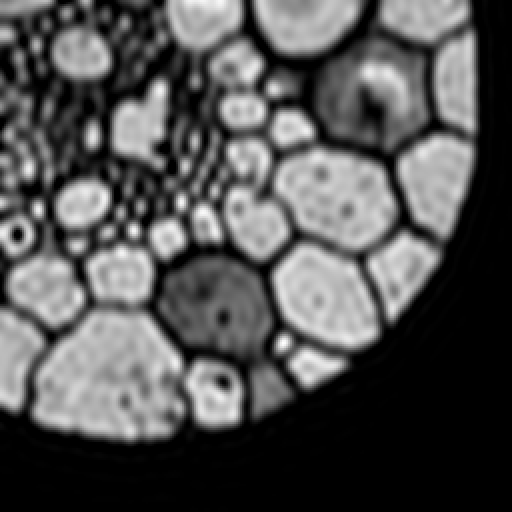}\\(d)
		\end{minipage}\\
		\begin{minipage}{0.24\textwidth}
			\centering
			\includegraphics[width = \textwidth]{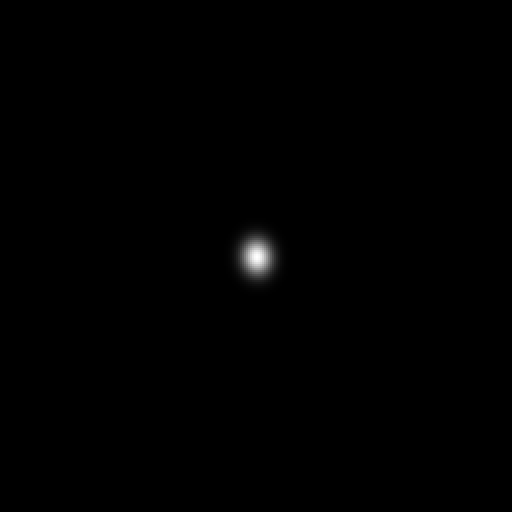}\\(e)
		\end{minipage}
		\begin{minipage}{0.24\textwidth}
			\centering
			\includegraphics[width = \textwidth]{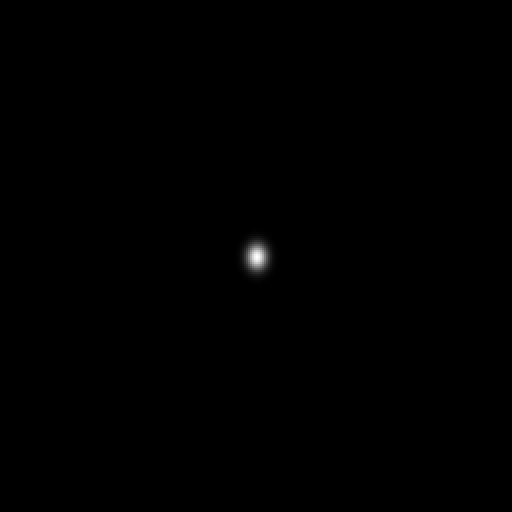}\\(f)
		\end{minipage}
		\begin{minipage}{0.24\textwidth}
			\centering
			\includegraphics[width=\textwidth]{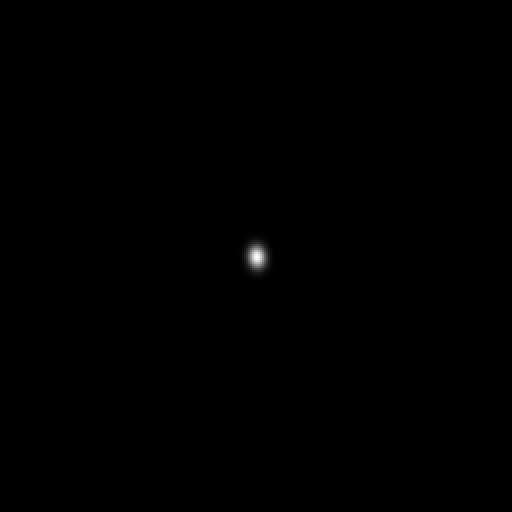}\\(g)
		\end{minipage}
		\begin{minipage}{0.24\textwidth}
			\centering
			\includegraphics[width = \textwidth]{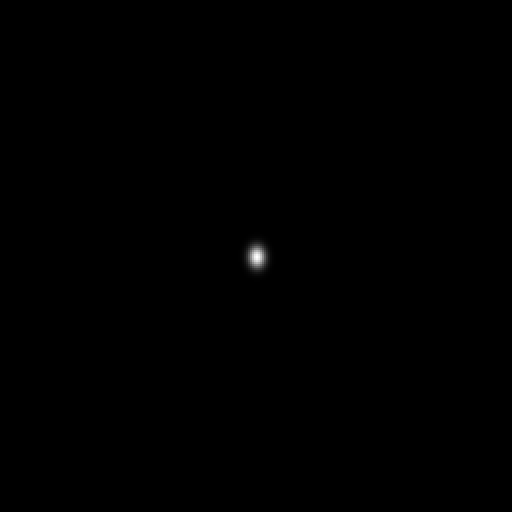}\\(h)
		\end{minipage}
		\caption{Reconstructions obtained by NLS$_{\ell_{2}-\ell_{p}}$ with $p = 1$. Panels (a-d) show reconstructions of the image $\bx$ and panels (e-h) reconstructions of the PSF at iterations 1,3, 5, and 9 of $\ell_p$ VarPro respectively from left to right.}
		\label{Figure: grain_rec}
	\end{figure}

	In this example we investigate the performance of NLS$_{\ell_2-\ell_p}$ in the reconstruction of the approximation of the forward operator as well as the desired image $\bx$. To do so, we consider the true image $\bx \in \R^{256\times 256}$ shown in Figure \eqref{fig:grain_problemSetup} (a). The image is then blurred with Gaussian PSF with parameters $\by_{\rm true} = (\sigma_1, \sigma_2, \rho) = (3, 4, 0.5)$ resulting in the PSF shown in Figure \eqref{fig:grain_problemSetup} (b), that after contamination with 1\% Gaussian noise result on the available measurements $\bd$ shown in Figure \eqref{fig:grain_problemSetup} (c). 
	We would like to determine an accurate approximation of this image from
	an available blur- and noise-contaminated version as well as an accurate approximation of the blurring operator. The actual parameters at each iteration are given in Table \ref{Tab: grainParams}.
	On this example we set $p = 1$ and the regularization matrix to be $\bL = \bW$, hence we aim to reconstruct solutions whose coefficients in the framelet domain $\bW$ are sparse. The reconstructed images and the reconstructed PSFs at iterations $i = 1,3,5, 9$ of $\ell_p$ VarPro are shown in Figure \eqref{Figure: grain_rec} panels(a-d) and (e-h) respectively. We conclude the numerical examples with the reconstructed parameters $\by = (\sigma_1, \sigma_2, \rho)$ at each of nine iterations, as well as the initial guess and $\by_{\rm true}$ that are shown in Table \eqref{Tab: grainParams}. They support the argument that we not only reconstruct good quality desired images, but good approximations for the desired parameters as well.
	
	\begin{table}[h!]
		\caption{Convergence of the parameters $\by = (\sigma_1, \sigma_2, \rho)$ at iteration $i$.}
		\label{Tab: grainParams}
		\begin{tabular}{|l|llllllllll|l|}\hline
			i               &  0      & 1      & 2      & 3      & 4      & 5      & 6      & 7      & 8      & 9    &$\by_{\rm true}$  \\\hline
			$\sigma_1$ & 5.0 & 4.15 & 3.65 & 3.37 & 3.17 & 3.04 & 3.04 & 3.04 & 3.04 & 3.03 &3.0\\
			$\sigma_2$ & 6.0 & 5.06 & 4.54 & 4.27 & 4.06 & 3.93 & 3.93 & 3.92 & 3.92 & 3.92 &4.0\\
			$\rho$      & 1.0 & 0.78 & 0.71 & 0.63 & 0.55 & 0.50 & 0.50 & 0.50 & 0.49 & 0.49 &0.5\\\hline
		\end{tabular}
	\end{table}
	
	\section{Conclusions and future work}\label{sec: conclusions}
	In this paper we propose a new iterative method to solve large-scale separable nonlinear inverse problems by combining a majorization-minimization technique in a generalized Krylov subspace method with the VarPro method in order to provide solutions with sparse and edge-preserving properties. The first contribution is in considering a general regularization term $\|\bL\bx\|_2^2$ and formulating its corresponding VarPro method. Further, we generalize it to $\|\bL\bx\|_p^p$, $0<p\leq 2$. Since the regularization parameter is known to play an important role and the nonlinear problem is sensitive to its choice, we adaptively tune the regularization parameter for every iteration by GCV. Furthermore, we carefully derive the Jacobian matrices for all cases considered and we motivate our choices with numerical considerations. Several numerical examples illustrated the performance of the proposed methods in terms of the quality of the reconstructed solution and the speed of the convergence. An advantage of the current approach is that it improves the quality of the reconstructed image $\bx$ as well as it provides a nice approximation of the parameters $\by$ that define the forward operator. Since we can obtain a more accurate solution when minimizing with respect to $\bx$, the typical behaviour of the nonlinear problem is that it will converge in fewer iterations.
	Further work will include: (a) more theoretical investigation of the proposed algorithms, (b) imposing nonnegativity constraint for the edge-preserving and sparsity promoting methods, (c) improvements of the computational cost of the Jacobian. Future potential applications of interest include superresolution imaging, nuclear magnetic resonance data analysis, machine learning, and nonlinear neural networks training, to mention a few. 
	
	\section{Acknowledgments}
	The authors would like to thank Julianne Chung for comments that helped improve the manuscript and for providing the code associated with the paper \cite{chung2010efficient}. 
	
	\bibliographystyle{plain}
	\bibliography{arxivref}
\end{document}